\documentclass[10pt]{article}

\begin{document}

\input amssym.tex

\newcommand{\Lor}{L^{\uparrow}_{+}}

\newtheorem{lem}{Lemma}
\newtheorem{defin}{Definition}
\newtheorem{theor}{Theorem}
\newtheorem{rem}{Remark}
\newtheorem{prop}{Proposition}
\newtheorem{cor}{Corollary}
\newenvironment{demo}
{\bgroup\par\smallskip\noindent{\it Proof: }}{\rule{0.5em}{0.5em}
\egroup}

\title{{Elements of Linear Algebra}\\
{\large	Lecture notes}}
\author{Ion I. Cot\u aescu\\ {\it The West University of Timi\c soara,}\\{\it V. 
		Parvan Ave. 4, RO-1900 Timi\c soara}}

\maketitle

\begin{abstract}
These pedagogical lecture notes address to the students in theoretical physics for helping them to understand the mechanisms of the linear operators defined on finite-dimensional vector spaces equipped with definite or indefinite inner products. The importance of the Dirac conjugation is pointed out presenting its general theory and a version of the Riesz theorem for the indefinite inner product spaces, based on the Dirac-Riesz map that combines the action of the Riesz map with that of the metric operator.  The matrix representations of the linear operators on vector spaces with definite or indefinite inner products is also presented.      
\end{abstract}

\newpage

\tableofcontents
\newpage

\section{Introduction}

The linear spaces are the most popular mathematical objects used in physics where the superposition principle is somewhat universal (apart from the non-linear classical systems). Thus the quantum physics works mainly with linear structures since  the state space is linear while the quantum observables are linear operators that form an operator algebra.    
These linear structures must offer us  quantities with a precise physical meaning that can be measured in various experiments. For this reason the mentioned linear structures must be equipped with norms or inner products able to generate real or complex numbers. Thus the principal linear structures we use in physics are the inner product real or complex vector spaces that can have a finite or even infinite number of dimensions. In what follows we restrict ourselves to the finite-dimensional ones.

On the other hand, the physical quantities must be in some sense invariant under the actions of  different symmetries that govern the systems under consideration. This means that the vector spaces must be seen as spaces carrying representations of the corresponding symmetry groups. It is known that in the case of the compact Lie groups there are unitary finite-dimensional representations in vector spaces with  positively defined inner products we refer here as unitary spaces. These are specific to the non-relativistic physics where the physical space is Euclidean and the quantum state space is a Hilbert one. 

However, in special or general relativity there are non-compact symmetry groups (as the Lorentz group of special relativity) whose finite-dimensional representations may be realized only in vector spaces with indefinite metric (or indefinite inner product).      
Moreover, the relativistic quantum physics in flat or curved backgrounds, in spite of the fact that  is based on similar principles as the non-relativistic one, introduces state spaces of the Krein type, having  indefinite inner products (called relativistic scalar products). This is the consequence of the existence of particles and antiparticle quantum modes in the entire relativistic quantum physics. Thus we see that the relativistic physics comply with various linear structures involving vector spaces with indefinite metric. 

For this reason, the present lecture notes are devoted to the general theory of the indefinite   
inner product spaces pointing out the role of the Dirac conjugation which is a general method able to offer suitable invariant quantities and self-adjoint operators when one consider non-unitary representations.  In addition, we present a version of the Riesz theorem  in indefinite inner product spaces based on a mapping that combine the Riesz map with the action of the metric operator. This map is called here the Dirac-Riesz map showing that this plays the similar role as the genuine Riesz one in the unitary case.  A special attention is paid to the matrix representations of the linear operators which is important in many applications in physics.

This paper is mainly of pedagogical interest helping the students in physics to become familiar with the usual linear algebra \cite{A1,A2,A3,A4}, or preparing them for studying advanced topics in linear algebra \cite{adA1,adA2}, group theory \cite{G1,G2,G3} or differential geometry \cite{GRG1,GRG2,GRG3} that represent the mathematical framework of the  actual theoretical physics.

\newpage

\section{Preliminaries}

The mathematical objects used in physics are quit complicated combining algebraic and topological properties defined at the level of the fundamental mathematics which remains far from the current needs of physical applications. For this reason, we start with a brief review of some important definitions that seems to be less familiar even for theoretical physicists.  In the same time we introduce the terminology and the principal conventions and notations.     

\subsection{Sets and mappings}

We denote by $A,\, B,\,...X,...$ different sets formed by elements $a,\,b,\,...x,\, y,...$. Two sets are equal, $A=B$, if have the same elements, i. e. $a\in A\,\Leftrightarrow a\in B$.  Given a set $A$ then an arbitrary element can be in $A$ ($a\in A$) or not in $A$  ($a\notin A$).  If every $a\in A$ is in another set, $a\in B$, the we say that $A$ is subset of $B$ and denote $A\subset B$ or $A\subseteq B$ when the possibility of $A=B$ is not excluded.  

The union of two sets is the set of elements that are in $A$ or in $B$. This is denoted as $A\cup B=\{x|x\in A \,{\rm or}\, x\in B,\,{\rm or\, both}\}$. Similarly, one defines the intersection of two sets, $A\cap B=\{x|x\in A \,{\rm and}\, x\in B\}$, and the difference $A\setminus B=\{x\in A|x\notin B\}$. We denote by $\emptyset$ the empty set (which has no elements) and we say that two sets, $A$ and $B$,  are disjoint if $A\cap B=\emptyset$. 

The product or Cartesian product of two sets is formed by ordered pairs $A\times B=\{(a,b)|a\in A,\, b\in B\}$. This operation is associative such that we can write $A\times (B\times C)=A\times B\times C$. The best-known example is of the Cartesian coordinates of a $n$-dimensional space that form ordered sequences of real numbers $(x^1,x^2,...,x^n)$ of the set ${\Bbb R}^n={\Bbb R}\times {\Bbb R}...\times {\Bbb R}$ where ${\Bbb R}$ is the set of real numbers. 

Given two sets $X$ and $Y$, an association $x\to y$  so that each  $x\in X$ is associated with only one element $y\in Y$ is called a function, mapping or map  being denoted by $f:X\to Y$ or  $X^{_{ \overrightarrow{~~~~~~} }^{~~f}}Y$. One says that $X$ is the domain of $f$ and ${\rm Im }f=f(X)\subset Y$ is the image of the map $f$. The image of a subset $U\subset X$ is denoted by $f(U)$. Other associations for which the domain is a subset of $X$ are called projections. 

A map $f:X\to Y$  can be:

1) injective when $x_1\not=x_2$ implies $f(x_1)\not=f(x_2)$;

2) surjective if for any $y\in Y$ there is at least one $x\in X$ such that $y=f(x)$; 

3) bijective when this is injective and surjective.

The composition of two maps, $f:X\to U$ and $g:U\to Y$  is the new map $g\circ f: X\to Y$ defined as  $(g\circ f) (x)=g(f(x))\in Y$ for all  $x\in X$ and $f(x)\in U$. 
This rule is illustrated  by the diagram 
\begin{equation}
\begin{array}{ccccc}
&f&&g&\\
X&\overrightarrow{~~~~~~~~~}&U&\overrightarrow{~~~~~~~~~}&Y
\end{array} =
\begin{array}{ccc}
&g\circ f&\\
X&\overrightarrow{~~~~~~~~~}&Y\,.
\end{array}
\end{equation}
The bijective map $f:X\to X$ is invertible and the inverse map $f^{-1}$ satisfies $f\circ f^{-1}=f^{-1}\circ f={\rm id}_X$ where  ${\rm id}_X$ is the identical map on $X$ acting as ${\rm id}_X(x)=x$ for all $x\in X$. For any invertible maps, $f$ and $g$,  we have $(f\circ g)^{-1}=g^{-1}\circ f^{-1}$ as represented in the next diagram, 
\begin{equation}\label{inverse}
\begin{array}{ccccc}
&f^{-1}&&g^{-1}&\\
X&\overleftarrow{~~~~~~~~~}&U&\overleftarrow{~~~~~~~~~}&Y
\end{array} =
\begin{array}{ccc}
&f^{-1}\circ g^{-1}&\\
X&\overleftarrow{~~~~~~~~~}&Y\,.
\end{array}
\end{equation}
When one composes different types of mappings  the following cases are possible\\ 
\begin{tabular}{rccccc}
~~~1)&any&$\circ$&injective&=&injective\\
2)&surjective&$\circ$ &any&=&surjective\\
3)&surjective&$\circ$ &injective&=&bijective
\end{tabular}\\
where 'any' means here one of  injective, surjective or bijective. 

The  composition of maps is associative, $f\circ(g\circ h)=(f\circ g)\circ h$.

\subsection{Algebraic structures}

We briefly review the principal algebraic structures considering {\em binary} operations (denoted by $+,\,\times,\,*,...$ etc.) which can be:

1)  {\em internal},  e. g.     $*:X\times X\to X$ such that $x*y\in X$ for all $x,y \in X$, 

2) {\em external},  e. g.  $\,\cdot :X\times Y\to X$ where $Y$ is another algebraic structure. Then $\alpha \cdot x \in X$ for all $x\in X$ and $\alpha \in Y$.

In general, a set can be equipped with one or more binary operations and some external ones but all these operations must be compatible among themselves according to a closed system of axioms. For this reason there are restrictions that select only few interesting algebraic structures that are of two types: the algebraic structures without an external operation and the linear structures whose external operation is the scalar multiplication.

The principal algebraic structures without this external operation are the groups, rings and numeric fields.
\begin{defin}{\bf: Group} is a set $X$ with only one internal operation, say $*$, which:

1) is associative, $x*(y*z)=(x*y)*z$
 
2)  has the {\em identity} element $e\in X$ so that $x*e=e*x=x$ for all $x\in X$,
  
3) for each $x\in X$ there is the inverse element $x^{-1}\in X$ which satisfies $x*x^{-1}=x^{-1}*x=e$.
\end{defin}
If, in addition, the operation is commutative, $x*y=y*x$ for all $x,y \in X$, then the group is called Abelian (or commutative). In this case the operation is addition, $+$, with the addition identity  denoted by $0\in X$ and the inverse of $x$ denoted by $-x$.\\   
\begin{defin}{\bf: Ring (with identity)} is an Abelian group under addition, $+$, equipped with multiplication $\cdot$ which:

1) is associative, $x\cdot(y\cdot z)=(x\cdot y)\cdot z$,

2) has a multiplicative identity $1$ such that $1\cdot x=x\cdot 1=x$,  

3) multiplication distributes over addition, $x\cdot (y+z)=x\cdot y + x\cdot z$ and $(x+y)\cdot z=x\cdot y + y\cdot z$.
\end{defin}
When the multiplication is commutative, $x\cdot y=y\cdot x$, we speak about  commutative rings.\\
\begin{defin}{\bf: Field} is a commutative ring (with identity), $X$, so that $X-\{0\}$ is an Abelian group under multiplication.
\end{defin}
Examples: the fields of rational (${\Bbb Q}$), real (${\Bbb R}$) or complex  (${\Bbb C}$) numbers. The complex numbers  are denoted as usually  by $z=\Re z+i\Im z$ where $\Re z,\,\Im z \in {\Bbb R}$. The complex conjugate of  $z$ is denoted as $\overline{z}=\Re z- i\Im z $ so that 
$|z|^2 = z\overline{z}=(\Re z)^2+(\Im z)^2$ where $|z|$ is the modulus of $z$. The same notation, $|x|$, stands for the absolute value of the real numbers $x\in{\Bbb R}$.

The linear structures are genuine vector spaces or vector spaces with additional binary operations.

\begin{defin}{\bf: Vector space}, $V$, over the numeric field ${\Bbb K}$  which  is  Abelian group under addition,  $+$, and closed with respect to the external operation of scalar multiplication, $\alpha x\in V$, $x\in V$,   $\alpha\in {\Bbb K}$,  that obeys the axioms:

1)  is distributive with respect to vector addition, $\alpha (x+y)=\alpha x+\alpha y$, 

2)  is distributive with respect to field addition,$(\alpha+\beta)x=\alpha x = \beta x$, 

3) is compatible with the field multiplication, $\alpha(\beta x)=(\alpha \beta) x$, 

4) there exists the identity element, $1 x=x$. 
\end{defin}
The null element of $V$ will be denoted by $0\in V$ in order to do not be confused with $0 \in {\Bbb K}$.\\ 
\begin{defin}{\bf: Algebra}, $A$, over  the numeric field ${\Bbb K}$ is a vector space over ${\Bbb K}$ equipped with an additional binary operation, $*$, which

1) is distributive to the left, $(x+y)*z=x*z+y*z$,

2) is distributive to the right $z*(x+y)=z*x+z*y$,

3) is compatible with the scalar multiplication, $(\alpha x)*(\beta y)=(\alpha \beta) (x*y)$.
\end{defin}
When the multiplication is associative, $x*(y*z)=(x*y)*z$, then the algebra is said associative.
The multiplication may be {\em symmetric} or {\em commutative} when $x*y=y*x$ or {\em skew-symmetric} if $x*y=-y*x$. If there exists the identity element $1$ such that $1*x=x*1=x$ for all $x\in A$ then we say that the algebra has unity.  The algebras with skew-symmetric multiplication can not have unity since then the algebra might be trivial $A=\{0\}$ (as long as  $x=-x \Rightarrow x=0$).

The principal virtue of the linear structures, vector spaces and algebras, is that the external operation of scalar multiplication allows one to define the bases as maximal systems of  linearly independent vectors.  
\begin{defin}\label{baseV} 
Suppose $V$ is a vector space over the field ${\Bbb K}$. An ordered system of vectors $e_1,e_2,...e_n \in V$ forms a basis if

1) the vectors  $e_1,e_2,...e_n$ are linearly independent,

2) any vector $x\in V$ is linearly expanded as 
\begin{equation}
x=\sum _{i=1}^n x^i e_i = x^i e_i \quad {\rm (Einstein's\,\, convention)} 
\end{equation}
where the numbers $x^i\in {\Bbb K}$ are the components (or coordinates) of the vector $x$ in this basis.
\end{defin}  
The (finite) number of the basis vectors defines the dimension of the vector space,  ${\rm dim} V=n$.  
\begin{rem}
We use the Einstein convention of dummy indices only for pairs of identical upper and lower indices. When their positions are ambiguous we use explicit sums.
\end{rem} 

The major challenge  is to find rational generalizations of this definition for $n\to \infty$ when we desire to have a {\em countable} basis $e_1,e_2,...e_i,..$ in which a vector $x$ should have a countable number of components $x^1,x^2,...x^i,..$.  Then we can calculate only finite linear combinations,    
\begin{equation}
x_{(n)}=\sum_{i=1}^n x^i e_i \quad n\in {\Bbb N}
\end{equation}
that form the sequence $\{x_{(n)}|n\in {\Bbb N}\}$. The problem is how could we define the convergence of this sequence, $x_{(n)}\to x$,  to the vector $x$ represented by the mentioned components.   Obviously, this cannot be done using algebra - in addition, we need topology.    

\subsection{Topologies} 

The minimal conditions in which the convergence of a sequence can be defined are offered by the topological spaces where the topology is given by open sets.
\begin{defin}
We say that $X$ is   topological space if this has a family of subsets $\tau(X)$, called open subsets, such that:

1) $X$ and and $\emptyset$ are open subsets,

2) any union of open sets is an open subset, 

3) the finite intersections of open sets are open subsets
\end{defin}
Then there are collections of open subsets $\{U_{\alpha}\}_{\alpha\in A}\subset \tau(X)$, labelled by some sets the indices, $A$, that form {\em covers} of $X$ obeying $X\subseteq \bigcup_{\alpha\in A}U_{\alpha}$.  Every point $p \in X$ is contained in some open set $U_p\in \tau(X)$ called neighbourhood of $p$. 

A boundary point of a set $S$ in a topological space $X$ is a point $p$ such 
that any of its neighbourhood $U_p$ contains both some points in $S$ 
and some points not in $S$.  The boundary $\partial S$ of the set $S$ is the set of boundary points of $S$.  When $\partial S\subset S$ then $S$ is closed.  The closed set $\overline{S}=S\cup \partial S$ is  called the closure of the set $S$. 
\begin{defin}
We say that the set $S\subset X$ is dense in $X$ if $\overline{S}=X$. 
\end{defin}
In general, a subset  may be neither open nor closed but there are sets that are considered closed and open simultaneously as  the space $X$ and the empty set.  One can show that if $A$ is an open set then its compliment $A^c=X\setminus A$ is closed.  

This framework is enough for defining the convergence of  sequences and the continuity of functions and algebraic operations. Let us start with the convergence.
\begin{defin}
Let $(x_1,x_2,... x_n,...)\subset X$ a sequence in the topological space $X$. We say that $x$ is the limit of this sequence, $x=lim_n x_n$,  if any  neighbourhood $U_x\in \tau(X)$ leaves outside a finite number of terms of this sequence. If the sequence has a limit then this is convergent in the topology of $\tau(X)$.
\end{defin}
Now  the continuity can be defined either using convergent sequences or staring with the independent definitions presented below. 
\begin{defin}
A map $f:X\to Y$ between two topological spaces is called continuous at $x\in X$ if for any neighbourhood $U_{f(x)}\in \tau(Y)$ of $f(x)$ there exists a neighbourhood $U_x \in \tau(X)$ so that $f(U_x)\subset U_{f(x)}$. A map is said continuous in a domain $D\subset X$ if this  is continuous in every point $x\in D$. When $D=X$ then the map is continuous.
\end{defin}
When the topological space is equipped with algebraic operations then these must be continuous with respect to the topology under consideration.
\begin{defin}
Let $X$ be a topological space equipped with the internal operation $*:X\times X\to X$. This operation is continuous in $X$ if for any $x,\,y\in X$ and any neighbourhood $U_{x*y}$ of $x*y$ there exists the neighbourhoods $U_x$ of $x$ and $U_y$ of $y$ such that for any $x'\in U_x$ and $y'\in U_y$ we have $x'*y'\in U_{x*y}$. 
\end{defin}
For the external operations the continuity is defined similarly.

The topological spaces are important from the mathematical point of view for establishing the most general conditions in which more complicated structures can be build using minimal sets of axioms. In physics we are interested only by spaces in which some general  topological properties appear as being natural. For this reason we briefly mention the definitions of the paracompact Hausdorff spaces which give rise to manifolds.   
\begin{defin}
If from any cover of open sets one can extract a finite one then the topological space is compact.
\end{defin}
\begin{defin}
$X$ is said  paracompact   if from any cover one can extract a locally finite one, in which every point of $X$ has a neighbourhood that intersects only finitely other sets in the cover. 
\end{defin}
Any compact space is paracompact.
\begin{defin}
When any two distinct points $p,p'\in X$ can be separated by neighbourhoods, $U_p$ and $U_{p'}$, so that  $U_p \bigcap U_{p'}=\emptyset$ then $X$ is a Hausdorff topological space.
\end{defin}
There are famous examples of spaces that are nor Hausdorff and neither paracompact but these are only of mathematical interest since in physics we consider only spaces having both these properties.

The typical paracompact  Hausdorff topological spaces are the metric spaces.
\begin{defin}
The metric topological space is the pair $(X,d)$ where $d: X\times X\to {\Bbb R}$ is the distance
function which satisfies

1) $d(x,y)=d(y,x)$, $\forall\, x,\,y\,\in X$  (symmetry),

2) $d(x,y)=0$ implies $x=y$,

3) $d(x,y)\le d(x,z)+d(z,y)$.  $\forall\, x,\,y,\,z\,\in X$ (triangle inequality)
\end{defin}
In a metric space $(X,d)$ one can introduce the topology of open balls $B(x,r)=\{y| d(x-y)<r\}$ with the center $x$ and any radius $r\in {\Bbb R}$. This topology is denoted  $\tau_d(X)$. Then the definitions of the convergence and continuity can be "translated" in terms of open balls. Thus we say that $x$ is the limit of the sequence $(x_1,x_2,... x_n,...)$ if for any $\epsilon >0$ there exists a finite number  $N(\epsilon)$ so that $x_n\in B(x,\epsilon)$ for all $n>N(\epsilon)$. The continuity can be reformulated in the same manner.

However, the metric spaces are important since therein on can define the basic notion of Cauchy sequence which can not be introduced without distance functions.
\begin{defin}
We say that the sequence $(x_1,x_2,... x_n,...)$ is a Cauchy sequence if for any $\epsilon>0$ there is $N(\epsilon)$ such that for $n,\,m>N(\epsilon)$ we have $d(x_n, x_m)<\epsilon$.
\end{defin}   
With this notion one can control the convergence of the sequences in infinite-dimensional spaces. 
\begin{defin}
A metric space is said complete if any Cauchy sequence is convergent.
\end{defin}
Any finite-dimensional metric space is anyway complete. The problems appear for the infinite-dimensional spaces where the completeness is a crucial axiom. 

The  further developments look for the context in which the distance functions can be defined on a natural way. For this purpose the space must be (at lest) a vector space with a topology given by a norm or by an inner product generating this norm.   Thus the principal topological spaces we meet in physics comply with the following hierarchy: topological vector space $>$ normed space $>$  inner product space. The infinite-dimensional complete normed spaces are the Banach spaces while the complete inner product spaces are called Hilbert spaces.

\section{Some linear algebra} 

The linear algebra is largely used in physics for describing states obeying the superposition principle, various classical or quantum physical observables and the geometric properties of the spacetime.  In what follows we review the principal  results of linear algebra we need for understanding the theory of classical and quantum fields on curved manifolds.  

\subsection{Duality and tensors}

The duality is a crucial and somewhat natural mechanism  that enables us to define correctly  different linear structures introduced in physics sometimes  in a intuitive but superficial manner. Our purpose here is to present the tensors as well-defined mathematical objects getting over the naive definition of  "quantities carrying indices which transform as...". 

\subsubsection{Linear mappings}

We consider vector space, denoted by $V$, $W$, $U...$, defined over a numeric field ${\Bbb K}$ (of  real,  ${\Bbb R}$ or complex, ${\Bbb C}$, numbers) and denote by $x,y,...\in V$ the vectors and by $\alpha, \beta...\in {\Bbb K}$ the scalars. 

Any subset of $V$  which is a vector space is called a subspace of $V$. The set $\{0\}$ and $V$ are trivial subspaces of $V$. Two subspaces $V_1$ and $V_2$ are said to be complementary subspaces  if $V_1\cup V_2=V$ and $V_1\cap V_2=\{0\}$. The dimensions of the non-trivial subspaces $V_1\subset V$ satisfy $1\leq{\rm dim}(V_1)<{\rm dim}(V)$ while the dimension of the complementary subspaces satisfy ${\rm dim} V_1+{\rm dim} V_2={\rm dim} V$. Moreover, any basis $e_1,e_2,...e_n$ in $V$ can be split in two disjoint bases  in $V_1$ and $V_2$ respectively. This property is somewhat similar to the orthogonal sum of subspaces but now we cannot speak about orthogonality since this has to be defined with respect to a inner product.

In general, the linear mappings can be defined between arbitrary  vector spaces  over the same field ${\Bbb K}$.
\begin{defin}
A map $f:V\to W$ is linear if 
\begin{eqnarray}
&&f(x+y)=f(x)+f(y) \quad \forall\, x,y\in V\nonumber\\
&&f(\alpha x)=\alpha f(x) \quad \forall\,\alpha\in {\Bbb K}\,,\,  x\in V\nonumber
\end{eqnarray}
\end{defin}
The identical map ${\rm id}_V:V\to V$  is linear,  the composition $f\circ g$ of  two linear maps, $f$ and $g$, is a linear map too and the inverse $f^{-1}$ of a bijective linear map $f$  is linear.   
\begin{lem}\label{linmap}
The linear maps $f:V\to W$ form a vector space, denoted by ${\rm Lin} (V,W)$, with respect to the addition and scalar multiplication of mappings,
\begin{eqnarray}
(f+f')(x)&=&f(x)+f'(x) \quad \forall\, x,\in V\nonumber\\
(\alpha f)(x)&=&\alpha f(x) \quad \forall\,\alpha\in {\Bbb K}\,,\,  x\in V\nonumber
\end{eqnarray}
\end{lem}
\begin{demo}
The null element is the null map $f_0$ giving $f_0(x)=0\in W$ for all $x\in V$ while the inverse of $f$ is $-f$ since $(f+(-f))(x)=f(x)+(-f(x))=0=f_0(x)$.      
\end{demo}

Each linear map $f:V\to W$ is related to a pair of important  vector subspaces: the kernel 
\begin{equation}
{\rm Ker} f =\{x\in V|f(x)=0\in W\}\subset V
\end{equation}
and the image  
\begin{equation}
{\rm Im} f =f(V)=\{y\in W| y=f(x),\, \forall x\in V\}\subset W\,.
\end{equation}
These are vector subspaces since $f$ is linear. In general ${\rm Im } f$ is a non-trivial subspace of $W$ since only the null map has ${\rm Ker } f_0=V$ and ${\rm Im } f_0=\{0\}$.  
It is shown that $f$ is surjective only if ${\rm Im} f=W$  and is injective when ${\rm Ker} f=\{0\}$.
A bijective map $f:V\to W$ has simultaneously  ${\rm Ker} f=\{0\}$ and  ${\rm Im} f=W$.
\begin{defin}
Two vector spaces, $V$ and $W$, are isomorphic if there exists a bijective linear map $f:V\to W$ 
\end{defin}
The isomorphic spaces must have the same dimension, ${\rm dim}(V)={\rm dim}(W)$. Reversely, any two vector spaces having the same (finite) dimension can be related among themselves through a bijective linear map. For example, giving two bases, $e_1,e_2,...e_n$ in $V$ and  $e'_1,e'_2,...e'_n$ in $W$,  there is a trivial isomorphism $\phi: V\to W$ so that $\phi(e_i)=e'_i$.  
\begin{defin}
Te linear maps $f:V\to V$ are called the endomorphisms of the space $V$  or simply,  linear operators.  The set of endomorphisms is denoted by ${\rm End}(V)$ 
\end{defin}
The  endomorphisms can be organized as a vector space over the field ${\Bbb K}$.  
The principal property of this space is that for any $f,h\in {\rm End}(V)$  the compositions $f\circ h$ and $h\circ f$ are also elements of ${\rm End}(V)$. Thus  ${\rm End}(V)$ is a vector space over ${\Bbb K}$ with the additional operation $\circ$ having  the following properties  
\begin{eqnarray}
(\alpha f+\beta g)\circ h&=&\alpha f\circ h +\beta g\circ h\\
h\circ (\alpha f+\beta g)&=&\alpha h\circ f +\beta h\circ g
\end{eqnarray}
where $f,g,h \in {\rm End}(V)$ and $\alpha, \,\beta\in {\Bbb K}$. With these operations, the vector space  ${\rm End}(V)$ takes a structure of an associative algebra but which is not commutative. In addition this has the unit mapping ${\rm id}_V$ such that ${\rm id}_V(x)=x$ for all $x\in V$. 
Thus the  conclusion is
\begin{rem}
The set  ${\rm End}(V)$ is an associative algebra with unity over the numeric field ${\Bbb K}$. 
\end{rem}
In this algebra, beside the bijective operators with ${\rm Ker} f=\{0\}$, we also find {\em degenerate} (or {\em singular}) operators $f$ having ${\rm Ker } f \not= \{0\}$  (with ${\rm dim}({\rm Ker} f)>0$ and ${\rm rank} f<{\rm dim} (V)$), including the null map $f_0$.     
\begin{defin}
The non-degenerate endomorphisms form the space ${\rm Aut}(V)$ of the {automorphisms} of the space $V$. 
\end{defin}
These mappings have the following obvious algebraic properties: if $f,g\in {\rm Aut}(V)$ then $f\circ g \in{\rm Aut}(V)$, for any $f\in{\rm Aut}(V)$ there exists the inverse $f^{-1}\in {\rm Aut}(V)$  and ${\rm id}_V\in {\rm Aut}(V)$.   
\begin{rem}
The set ${\rm Aut}(V)$ has the structure of a group with respect to the operation composition of maps.
\end{rem}
Under such circumstances the set ${\rm Aut}(V)$ is no longer a vector space since this cannot include degenerate operators.

\subsubsection{Duality}

A special case of linear maps is  when the image is just the field ${\Bbb K}$ considered as one-dimensional vector space defined over itself. 
\begin{defin} The linear map $\hat y: V\to {\Bbb K}$ is called a linear functional (or form).
\end{defin}
The linear functionals defined on $V$,  denoted from now by $\hat x,\hat y,...$, are called covectors. Any functional $\hat y:V\to {\Bbb K}$  maps a vector $x\in V$ into the number $\hat y(x)\in {\Bbb K}$. The covectors form the {\em dual} space $V^*={\rm Lin} (V,{\Bbb K})$ of $V$ which is a vector space over the same field  ${\Bbb K}$  with respect to the  addition and scalar multiplications, 
\begin{eqnarray}
(\hat y +\hat y')(x)&=&\hat y(x)+\hat y'(x) \quad \forall x\in V\\
(\alpha \hat y)(x)&=&\alpha \hat y(x)  
\end{eqnarray}
as it results from lemma  (\ref{linmap}). Reversely, the vectors  $x\in V$ can be seen as linear functionals  $x:V^*\to {\Bbb K}$ on $V^*$.
\begin{lem}
There exists a natural isomorphism between the spaces $V$ and ${\rm Lin} (V^*,{\Bbb K})$ called the canonical isomorphism.
\end{lem}
\begin{demo}
We define this isomorphism,  $\phi: V\to  {\rm Lin} (V^*,{\Bbb K})$, so that $(\phi(x))(\hat y)=\hat y(x)\in {\Bbb K}$
\end{demo} Thus we find the  rule $(V^*)^*\cong V$. 
\begin{defin}
The map  $\langle~,~\rangle : V^*\times V\to {\Bbb K}$ defined as $\langle \hat y,x\rangle=\hat y(x)$ for all $x\in V$ and $\hat y\in V^*$ is called the dual form of the pair of dual spaces $(V,V^*)$.
\end{defin}
The dual form is linear in both their arguments,
\begin{eqnarray}
\langle \alpha\hat y +\beta \hat y', x\rangle&=&\alpha \langle \hat y,x\rangle+\beta \langle \hat y',x\rangle\,, \\
\langle \hat y, \alpha x  +\beta  x'\rangle&=&\alpha \langle \hat y,x\rangle+\beta \langle \hat y,x'\rangle\,. 
\end{eqnarray} 
Notice that the dual form should not be confused with scalar or inner products which are defined on $V\times V$ and have supplementary symmetries.

\begin{defin}\label{fstar}
If $f:V\to W$ is a linear map then the map $f^*:W^*\to V^*$ which satisfies 
$\langle f^*(\hat y),x\rangle=\langle \hat y, f(x)\rangle$ for all   $x\in V$ and $\hat y \in W^*$ is called the dual map of $f$.
\end{defin}
It is no difficult to show that $f^*$ is also a linear map uniquely determined by the map $f$ so that the spaces ${\rm Lin} (V,W)$ and ${\rm Lin} (W^*,V^*)$ are isomorphic. We note that  if $f=f_0$ is the null map of ${\rm Lin} (V,W)$ then $f^*$ is the null map of ${\rm Lin} (W^*,V^*)$.   Similar properties have the pairs of isomorphic spaces ${\rm End}(V)$ and  ${\rm End}(V^*)$ and respectively  ${\rm Aut}(V)$ and  ${\rm Aut}(V^*)$.  Moreover, the following calculation rules are obvious: $(f^*)^*=f$ and
\begin{eqnarray}
&&(\alpha f +\beta g)^*=\alpha f^* +\beta g^*\,, \quad \alpha, \beta \in {\Bbb K}\,,\\
&&(f\circ g)^*=g^*\circ f^*\,.
\end{eqnarray}
\begin{defin}\label{selfdual}
The linear mappings $f\in {\rm Lin}(V.V^*)$ are called sefldual denoting $f^*=f$.
\end{defin}

Let us consider now a basis $e=\{e_1,\,e_2, ... e_n\}\subset V$   as in definition  (\ref{baseV})
and the set  $\hat e=\{\hat e^1,\,\hat e^2, ... \hat e^n\} \subset V^*$ that accomplishes the {\em duality} condition
\begin{equation}\label{dualb}
\hat e^i(e_j)=\langle\hat e^i, e_j\rangle=\delta_j^i\,.
\end{equation}
\begin{theor}
Given the basis $e \subset V$  then Eq. (\ref{dualb}) uniquely determines the set $\hat e \subset V^*$ of linearly independent covectors which forms a basis in $V^*$ called the dual basis of the basis $e$. 
\end{theor}
\begin{demo}
The system of covectors  is linearly independent since the equation $\alpha_i \hat e^i=0$ has only trivial solutions since $\alpha_i \hat e^i(e_j) =\alpha _j=0$ for $j=1,2,...n$.  For any $x=x^j e_j\in V$ and  $\hat y\in V^*$ we have $\hat y (x)=\hat y (e_i)x^i=\hat y_i x^i$ but the same result is obtained starting with the expansion $\hat y=\hat y_i \hat e^i$ since then $\hat y(x)=\hat y_i x^j \hat e^i(e_j)=\hat y_i x^i$.  Therefore, $\hat e$ is a basis in $V^*$
\end{demo}
\begin{rem}
In what follows we consider only  pairs (or systems) of dual bases, $(e,\hat e)$ in the spaces $V$ and $V^*$.  
\end{rem}
In a given system of dual bases  any vector $x\in V$ may be expanded as 
\begin{equation}\label{expv}
x=x^1e_1+x^2e_2+...\,  x^n e_n =  e_i \langle \hat e^i,x\rangle
\end{equation}
in terms of its  components (or coordinates) in this basis, $ x^i =\langle \hat e^i,x\rangle \in {\Bbb K}$. The covectors $\hat y\in V^*$ can also be expanded  in terms of their components (or coordinates) $\hat y_i = \langle \hat y,e_i\rangle \in {\Bbb K}$ in the basis $\hat e$, 
\begin{equation}\label{expcv}
\hat y=\hat y_1 \hat e^1+\hat y_2 \hat e^2+...\,\hat y_n \hat e^n = \langle \hat y,e_i\rangle \hat e^i\,.
\end{equation}
Then, according to Eq. (\ref{dualb}), we find the expansion of the dual form, 
\begin{equation}
\langle \hat y, x\rangle = \hat y_1 x^1 +\hat y_2 x^2 +...\,\hat y_n x^n  = \langle \hat y,e_i\rangle \langle \hat e^i,x\rangle
\end{equation}
in the system of dual bases $(e, \hat e)$. 

\subsubsection{Tensors}

The natural generalization of the linear maps are the multilinear maps defined on different Cartesian products of arbitrary vector spaces. Let $V_1,\, V_2,\,...V_n$ and $W$ be vector spaces over the field ${\Bbb K}$ and denote by $x_i,\in V_i$ the vectors of any $V_i$  ($i=1,2,...n$).
\begin{defin}\label{multil}
A  map $f:V_1\times V_2...\times V_n \to W$ is multilinear if this is linear in each of its arguments,
\begin{equation}
f(....,\alpha x_i+\beta x'_i,.....)=\alpha f(..., x_i,....)+\beta f(..., x'_i,...)\,. 
\end{equation} 
\end{defin}
The set of multilinear maps can be organized as a vector space over ${\Bbb K}$ denoted by ${\rm Lin}(V_1,V_2,..., V_n; W)$ the dimension of which is 
\begin{equation}
{\rm dim}\,{\rm Lin}(V_1,V_2,..., V_n; W)={\rm dim} V_1{\rm dim} V_2 ...{\rm dim} V_n{\rm dim} W.
\end{equation}
When $W={\Bbb K}$ we say that ${\rm Lin}(V_1,V_2,..., V_n; {\Bbb K})$ is the space of multilinear forms (or functionals) over $V_1\times V_2 ...\times V_n$. These mappings allow one to construct the different classes of tensors using a specific operation, $\otimes$,  called the tensor product.
\begin{defin}
The tensor product of a set of vector spaces $\{V_i|i=1,2...n\}$ is  the map  
\begin{equation}
\otimes: V_1\times V_2...\times V_n \to V_1\otimes V_2... \otimes V_n\cong  {\rm Lin}(V^*_1,V^*_2,...V_n^*; {\Bbb K})
\end{equation} 
such that for any $x_i\in V_i$, and $\hat y_i \in V_i^*$ we have 
\begin{equation}
(x_1\otimes x_2...\otimes x_n) (\hat y_1, \hat y_2,... \hat y _n)=\langle \hat y_1, x_1\rangle \langle \hat y_2, x_3\rangle ... \langle \hat y_n, x_n\rangle\,.
\end{equation}
\end{defin}
The tensor product  space $V_1\otimes V_2...\otimes V_n$ includes all the linear combinations of the tensor product vectors, $x_1\otimes x_2...\otimes x_n$, since the tensor product is a multilinear operation in the sense of definition (\ref{multil}). On the other hand,  the tensor product  may be seen as a collective operation where only the ordering of the vector spaces is significant.  This is to say that the tensor product is an {\em associative}  operation that can be performed in any order but keeping unchanged the ordering of the vector spaces.
\begin{defin}
Given the bases $e_{(i)}\subset V_i$,  the set  
\begin{equation}\label{tensbase}
\{e_{(1) i_1}\otimes  e_{(2) i_2}...\otimes  e_{(n) i_n}| i_k=1,2,...{\rm dim} V_k,   k=1,2,...n \}
\end{equation}
forms the tensor product basis of the space $V_1\otimes V_2...\otimes V_n$.  
\end{defin}
However, this is not the unique option since many other bases can be defined using linear combinations.  The dual bases can be constructed in the same manner since there are  almost obvious duality properties as given by the next Lemma. 
\begin{lem}\label{lemdual}
There exists the natural isomorphism 
\begin{equation}
(V_1\otimes V_2... \otimes V_n)^*\cong  V^*_1\otimes V^*_2... \otimes V^*_n \cong {\rm Lin}(V_1,V_2,...V_n; {\Bbb K})\,.
\end{equation}
\end{lem} 
Therefore, the dual basis of the base (\ref{tensbase}) that reads
\begin{equation}\label{dtensbase}
\{\hat e^{(1) i_1}\otimes \hat e^{(2) i_2}...\otimes \hat e^{(n) i_n}| i_k=1,2,...{\rm dim} V_k,   k=1,2,...n \}
\end{equation}
represents the tensor product basis of the space $(V_1\otimes V_2...\otimes V_n)^*$.  

An important particular example helps us to understand how different algebraic object are naturally related.
\begin{lem} 
For any two vector spaces $V$ and $W$ over the field ${\Bbb K}$ there exists the natural isomorphism 
$V\otimes W^*\cong {\rm Lin}(V^*,W;{\Bbb K})\cong {\rm Lin}(V,W)$.
\end{lem}
\begin{demo}
Let $\phi:  {\rm Lin}(V^*,W;{\Bbb K})\to {\rm Lin}(V,W)$ so that for any $\tau \in {\rm Lin}(V^*,W;{\Bbb K})$ the map $f_{\tau}=\phi(\tau)$ satisfies $\tau(\hat x, f_{\tau} (x) )=\langle \hat x,x\rangle$ for all $x\in V$ and $\hat x \in V^*$. Then it is a simple exercise to verify that $\phi$ is an isomorphism. 
\end{demo}

The tensor product spaces can be related among themselves through various multiliniar mappings but the good tools are the tensor products of mappings.
\begin{defin}
Let us suppose that $n$ linear maps $f_k:V_k\to W_k$, $k=1,2,...n$ are known.  Then we say that that multilinear map $f=f_1\otimes f_2...\otimes f_n$ is the tensor product map if 
\begin{equation}
f(x_1\otimes x_2...\otimes x_n)= f_1(x_1)\otimes f_2(x_2)...\otimes f_n(x_n)\quad \forall x_i\in V_i\,.
\end{equation} 
\end{defin}
Then we denote 
\begin{eqnarray}
 f&\in& {\rm Lin}(V_1\otimes V_2...\otimes V_n, W_1\otimes W_2...\otimes W_n)\nonumber\\
 &\cong&
 {\rm Lin}(V_1,W_1)\otimes {\rm Lin}(V_2,W_2)...\otimes {\rm Lin}(V_n,W_n)\,. 
\end{eqnarray}

Particular tensor products are those involving only copies of a given vector space $V$ and its dual $V^*$. These give all  the tensors associated to the vector space $V$.
\begin{defin}
The vector space of the contravariant tensors  of rank $n$ is defined as 
\begin{equation}
T^{n}(V)=V^{\otimes n}=\underbrace{V\otimes V...\otimes V}_n = \,.
\end{equation}
while the space of the covariant tensors of rank $n$ is 
\begin{equation}
T_{n}(V)=V^{*\,\otimes n}= \underbrace{V^*\otimes V^*...\otimes V^*}_n\,.
\end{equation}
\end{defin}
According to lemma (\ref{lemdual}) these spaces are dual each other, $[T^n(V)]^*=T_n(V)$. Notice that these spaces have the same dimension and therefore these are isomorphic through the trivial isomorphism of their dual bases but they must be seen as different mathematical objects since the covariant and contravariant tensors are different multilinear forms.  In particular, we  identify $T^1(V)=V$ and $T_1(V)=V^*$.

Giving a pair of dual bases $(e,\hat e)$ of $(V,V^*)$ we can construct the tensor product bases 
\begin{eqnarray}
e^{\otimes n}&=&\{ e_{i_1}\otimes  e_{i_2}...\otimes  e_{i_n}|i_k=1,2,...{\rm din} V\}\\
\hat e^{\otimes n}&=&\{ \hat e^{i_1}\otimes \hat e^{i_2}...\otimes \hat e^{i_n}|i_k=1,2,...{\rm din} V\}
\end{eqnarray}
of the dual spaces $T^n$ and $T_n$.  In general, the spaces of  tensors of arbitrary ranks, 
\begin{equation}
T^{n_1\, \cdot\, n_2...}_{\,\,\cdot\, m_1\cdot\, m_2...}(V)=T^{n_1}(V)\otimes T_{m_1}(V)\otimes T^{n_2}(V)\otimes T_{m_2}(V).... \,,
\end{equation} 
have  tensor product bases of the form $ e^{\otimes n_1}\otimes \hat e^{\otimes m_1}\otimes  e^{\otimes n_2}....$. 
For example, the tensors of the space  $T^{n}_{\,\cdot\, m}(V)=T^{n}(V)\otimes T_{m}(V)$  are the multiliniar forms 
\begin{equation}
\tau:\underbrace{V^*\times V^* ....\times V^*}_{n}\times \underbrace{V\times V ....\times V}_{m}\to {\Bbb K}\,.
\end{equation}
The tensor product basis  $ e^{\otimes n}\otimes \hat e^{\otimes m}$ of this space 
is formed now by tensors that  act on the basis vectors and covectors as 
\begin{equation}
(e_{i_1}\otimes e_{i_2}...\otimes e_{i_n}\otimes \hat e^{j_1}\otimes \hat e^{j_2}...\otimes \hat e^{j_m})(\hat e^{k_1}, \hat e^{k_2},..., e_{l_1},e_{l_2},....)=\delta^{k_1}_{i_1}...\delta^{j_1}_{l_1}...
\end{equation}
Any tensor $\tau \in T^{n}_{\,\cdot\, m}(V)$ can be expanded in this basis, 
\begin{equation}
\tau=\tau^{i_1i_2...i_n}_{\,\,\cdot\,\,\cdot\,\,\cdots\,\cdot\, j_1j_2...j_m}e_{i_1}\otimes e_{i_2}...\otimes e_{i_n}\otimes \hat e^{j_1}\otimes \hat e^{j_2}...\otimes \hat e^{j_m}
\end{equation}
laying out its components,
\begin{equation}
\tau^{i_1i_2...i_n}_{\,\,\cdot\,\,\cdot\,\,\cdots\,\cdot\, j_1j_2...j_m}=\tau(\hat e^{i_1}, \hat e^{i_2},..., e_{j_1},e_{j_2},....)\,.
\end{equation}

A special map can be easily defined starting with the tensor components. 
\begin{defin}
The map ${\rm Tr}^k_l :T^{n}_{\,\cdot\, m}(V)\to T^{n-1}_{\,\,\,\cdot\,\,\, m-1}(V)$ is a contraction of the $k$-th upper index with the $l$-th lower one if 
\begin{eqnarray}
&&{\rm Tr}^k_l (e_{i_1}\otimes e_{i_2}...\otimes e_{i_k}...\otimes e_{i_n}\otimes \hat e^{j_1}\otimes \hat e^{j_2}...\otimes \hat e^{j_l}...\otimes \hat e^{j_m})\nonumber\\
&&~~~~~~~~=\langle\hat e^{i_k}, e_{j_l}\rangle\, e_{i_1}\otimes e_{i_2}...\otimes e_{i_n}\otimes \hat e^{j_1}\otimes \hat e^{j_2}...\otimes \hat e^{j_m}\nonumber\\
&&~~~~~~~~=\delta^{i_k}_{j_l}\,e_{i_1}\otimes e_{i_2}...\otimes e_{i_n}\otimes \hat e^{j_2}\otimes \hat e^{j_1}...\otimes \hat e^{j_m}
\end{eqnarray}
\end{defin}
This means that the components of the tensor $\theta ={\rm Tr}^k_l (\tau)$ are sums over the mentioned indices of the components of the tensor $\tau$, 
\begin{equation}
\theta^{i_1...i_n}_{\,\,\,\cdots\,\cdots\, j_1...j_m}=\tau^{i_1...s...i_n}_{\,\,\cdot\,\,\cdots\,\,\cdot\,\,\cdots\, j_1...s...j_m}
\end{equation}   
Many contractions can be performed successively decreasing the tensor rank as we desire. For example the tensors $T^n_{\,\cdot\, n}(V)$ can be contracted into scalars by the {\em  trace} mapping   ${\rm Tr}: T^n_{\,\cdot\, n}(V)\to {\Bbb K}$ giving ${\rm Tr}(\tau)=\tau^{i_1...i_n}_{\,\,\cdot\,\,\cdots\,\,\cdot\,\, i_1....i_n}\in {\Bbb K}$.

\subsection{Matrix algebra}

In this section we briefly review the theory of matrices that are omnipresent in physics  where linear structures or linear transformations are used. Moreover, the arithmetic vector spaces formed by column or line matrices are the archetypes of all the vector spaces and the source of the Dirac bra-ket formalism.

\subsubsection{Matrices and determinants}

We denote by $ A ,\, B ,...$ any  $m\times n$ matrices, with $m$ lines and $n$ columns. The matrix elements of a $m\times n$ matrix $ A $  are denoted as $A^i_j$ where $i=1,2,...m$ is the line index and $j=1,2,...n$ is  the column one. The matrix elements are numbers from a field ${\Bbb K}$ that can be ${\Bbb K}={\Bbb R}$ when we have real matrices and ${\Bbb K}={\Bbb C}$ for complex matrices.

Two binary operations preserve the matrix dimension, i. e. the addition
\begin{equation}
 A + B = C \quad\Rightarrow\quad A^i_j+B^i_j=C^i_j
\end{equation}
and the scalar multiplication
\begin{equation}
\alpha A = B \quad\Rightarrow\quad \alpha A^i_j=B^i_j\,.
\end{equation}
The matrix multiplication,
\begin{equation}\label{multm}
 A \, B = C \quad\Rightarrow\quad A^i_kB^k_j=C^i_j
\end{equation}
changes the  dimension according to the well-known rule
\begin{equation} 
 (m\times
n')(n'\times n)=(m\times n). 
\end{equation}
Therefore, in Eq. (\ref{multm}) we must take $i=1,2,...m$, $j=1,2,...n$ while
the summation index $k$ ranges from $1$ to $n'$. The matrix multiplication is associative, 
$A(BC)=(AB)C$.

In the case of the complex  matrices ( ${\Bbb K}={\Bbb C}$) the Hermitian conjugation plays the same role as the simple transposition of the real matrices.   
\begin{defin}
The Hermitian conjugated $ A ^+$ of a complex $m\times n$ matrix $ A $  is a $n\times m$ matrix with the matrix elements $(A^{+})^{i}_j=\overline{A^j_i}$. For the real matrices this  reduces to the transposition $(A^{T})^{i}_j=A^j_i$.
\end{defin}
Then the following properties are obvious: $(A^+)^+=A$ and
\begin{equation}
(A+B)^+=A^++B^+\,,\quad (A B)^+=B^+ A^+\,,\quad (\alpha A)^+=\overline{\alpha} A^+\,, \quad \alpha \in {\Bbb C}\,.
\end{equation}
\begin{defin}
The complex $n\times n$ matrices which satisfy $A^+=A$ are called Hermitian or self-adjoint  matrices while those having the property $A^+=A^{-1}$ are unitary matrices.  
\end{defin}
The real matrices with a similar properties are the symmetric matrices, $A^T=A$, and the {\em orthogonal} ones, $A^T=A^{-1}$.   

An important invariant of a $n\times n$ matrix $A$ is its determinant which can be calculated according to the rule,
\begin{equation}
{\rm det}\, A =\varepsilon_{i_1 i_2 ... i_n}A^{i_1}_1A^{i_2}_2...A^{i_n}_n\,, 
\end{equation}
where $\varepsilon_{i_1 i_2 ... i_n}$ is the totally anti-symmetric Levi-Civita symbol (with $\varepsilon_{123...n}=1$ and changing the sign to each index transposition). Then one can prove that ${\rm det}\,{\bf 1}_n=1$, ${\rm det}(\alpha A)=\alpha^n {\rm det}\, A$ and  the next theorem. 
\begin{theor}\label{detAB}
Suppose $A,\,B$ be $n\times n$ matrices. Then  ${\rm det} (AB)= {\rm det}\, A\, {\rm det}\,B$.
\end{theor} 
The principal consequence is that the invertible $n\times n$ matrices must have non-zero determinants since ${\rm det}\, A^{-1}=({\rm det}\, A)^{-1}$.
\begin{defin}
A $n\times n$ matrix $A$ with ${\rm det}\, A=0$ is called a singular matrix. The invertible matrices are non-singular.
\end{defin}     

In the case of the $n\times n$ complex matrices we find ${\rm det}A^+=\overline{{\rm det} A}$ and for the real matrices we have ${\rm det}A^T={{\rm det} A}$. Therefore,  if $A$ is a unitary matrix then ${\rm det} A^+=\overline{{\rm det} A}= 
({\rm det} A)^{-1}$ and consequently $|{\rm det} A|^2=1$. Similarly, if $A$ is an orthogonal matrix then  $({\rm det} A)^2=1$ which means that ${\rm det} A=\pm 1$. We proved thus the next theorem.
\begin{theor}
The determinants of the unitary matrices are arbitrary phase factors while  the determinants of the orthogonal matrices take the values  $\pm 1$.
\end{theor}
This result is important since it will help us to classify the different subgroups of the group $GL(n, {\Bbb K})$. 

\subsubsection{Arithmetic vector spaces}\label{Aritspaces}

The simplest algebraic structure involving matrices are the arithmetic vector spaces we present by using the traditional "bra-ket" notation introduced by Dirac and largely used in quantum mechanics.
\begin{defin}
We say that the set of column matrices
$|x\rangle=(x^1,x^2,...x^n)^T$ with the components $x^i\in {\Bbb K}$ constitutes the $n$-dimensional arithmetic vector space over ${\Bbb K}$ denoted by ${\Bbb K}^n$
\end{defin}
The addition ($|z\rangle = |x\rangle +|y\rangle \Rightarrow z^i=x^i+y^i$) and the scalar multiplication 
($|y\rangle=\alpha|x\rangle \Rightarrow y^i=\alpha x^i$) reduce now to simple matrix operations.  
The vectors $|x\rangle \in {\Bbb K}^n$, called in physics {\em ket} vectors,  can be expanded  in the standard or {\em natural} basis of this space, known as the {ket} basis,
\begin{equation}
|1\rangle=\left(
\begin{array}{c}
1\\
0\\
\cdot\\
0
\end{array}\right)\,, \quad
|2\rangle =\left(
\begin{array}{c}
0\\
1\\
\cdot\\
0
\end{array}\right)\,,...\quad
|n\rangle =\left(
\begin{array}{c}
0\\
0\\
\cdot\\
1
\end{array}\right) \,.
\end{equation}
The dual vector space $({\Bbb K}^n)^*$ of ${\Bbb K}^n$ is formed by the covectors (or {\em bra} vectors) $\langle\hat y|=(\hat y_1,\hat y_2,...\hat y_n)$ and has the standard dual basis, or the  bra basis,   
\begin{eqnarray}
\langle 1|&=&(1,0,...,0)\,,\nonumber\\
\langle 2|&=&(0,1,...,0)\,,\\ 
&\vdots &\nonumber\\
\langle n|&=&(0,0,...,n)\,,\nonumber
\end{eqnarray} 
that satisfy the duality condition $(\langle i|)(|j\rangle)= \langle i|j\rangle=\delta^i_j$ (calculated using matrix multiplication) and have the natural property $(|i\rangle)^T=\langle i|$.  

Hereby we see that the vector components can be derived using the matrix multiplication  $x^i=(\langle i|)(|x\rangle)=\langle i|x\rangle$ and similarly for those  of the covectors  $\hat y_i=\langle \hat y|i \rangle $.  Notice that the indices of the ket basis are considered here as lower indices while those of the bra basis are seen as upper ones. Nevertheless, since the bra-ket notation does not point out the index position we write the explicit sums when the Einstein convention could lead to confusion. Thus the identity matrix ${\bf 1}_{n}\in {\rm End} ({\Bbb K}^n)$ can be put in the form 
\begin{equation}
{\bf 1}_n={\rm diag}(1,1,...1)=\sum_{i=1}^n |i\rangle\langle i|\,,
\end{equation} 
allowing us to write the  expansions of the vectors and covectors,  
\begin{equation}
|x\rangle ={\bf 1}_n |x\rangle =\sum_{i=1}^n  |i\rangle\langle i|x\rangle\,, \quad 
\langle \hat x| =\langle \hat x|{\bf 1}_n  =\sum_{i=1}^n\langle \hat x|i\rangle \langle i| \,,
\end{equation}
and that of the dual form  $\langle \hat y|x\rangle=\langle \hat y|{\bf 1}_n|x\rangle=\sum_i\langle \hat y|i\rangle\langle i |x\rangle$.

In general, any linear mapping between two arithmetic vector spaces is a matrix $M$  transforming ket vectors, $|x\rangle \to |x'\rangle=M|x\rangle$, while the natural bases remain fixed.  For this reason we do not use distinctive notations for the linear maps and their matrices when we work with arithmetic vector spaces. Thus we denote by $A\in {\rm Lin}({\Bbb K}^n, {\Bbb K}^m)$ the  linear map, $A:{\Bbb K}^n\to {\Bbb K}^m$, and its $m\times n$ matrix. The dual map $A^*:({\Bbb K}^m)^*\to ({\Bbb K}^n)^*$ has a different action but the same matrix, $A$. Therefore, we can use the notations
\begin{equation}\label{AAstar}
A(|x\rangle)=A|x\rangle\,, \quad A^*(\langle \hat y|) =\langle \hat y| A\,,
\end{equation}
since $\langle \hat y|A(|x\rangle)=A^*(\langle \hat y|) |x\rangle=\langle \hat y|A|x\rangle$ for all 
$|x\rangle \in {\Bbb K}^n$ and $\langle \hat y|\in ({\Bbb K}^m)^*$. We recover thus the well-known  {\em direct} right and {\em dual} left actions of the linear mappings in the bra-ket formalism.

Any linear equation, $|y\rangle =A|x\rangle$, can be expanded in  the natural bases, as $\langle i|y\rangle =\langle i|A {\bf 1}_n|x\rangle=\sum_k \langle i|A|k\rangle\langle k |x\rangle$, obtaining  the linear system $y^i= A^i_k x^k$ written with the matrix elements 
\begin{equation}
A^i_j=\langle i|A|j\rangle =\left(0,...1,... 0\right)\left(
\begin{array}{cccc}
A^1_1&A^1_2&\cdots &A^1_n\\
A^2_1&A^2_2&\cdots &A^2_n\\
\vdots &\vdots &A^i_j&\vdots\\
A^m_1&A^m_2&\cdots &A^m_n
\end{array}\right)
\left(
\begin{array}{c}
0\\
\vdots\\
1\\
\vdots
\end{array}\right) 
\end{equation}
In general, one may use the expansion   
\begin{equation}
A={\bf 1}_mA{\bf 1}_n=\sum_{i=1}^m\sum_{j=1}^n |i\rangle\langle i|A|j\rangle\langle j|\,. 
\end{equation}

The linear operators $A,B,...\in{\rm End} ({\Bbb K}^n)$  are now simple $n\times n$  matrices. 
\begin{lem}
The set ${\rm End} ({\Bbb K}^n)$ forms  an associative algebra with respect the addition, scalar multiplication and matrix multiplication.
\end{lem}
\begin{demo}
The $n\times n$ matrices form a vector space with respect to the matrix addition and scalar multiplication (with the zero $n\times n$ matrix $|0|$ having only zero matrix elements). The matrix multiplication is associative and there is the identity matrix  ${\bf 1}_n$
\end{demo} 

The singular operators  of ${\rm End} ({\Bbb K}^n)$ are {\em projection} operators. The simpler ones are the elementary projection operators, $P_k=|k\rangle\langle k|$ (no sum), which project any vector $|x\rangle \in {\Bbb K}^n$ to the one-dimensional subspace $V_k=\{\alpha |k\rangle | \alpha \in {\Bbb K}\}$ and  have the following properties
\begin{eqnarray}
1)\, {\rm orthogonality:}&\quad & P_j P_k =\left\{\begin{array}{ccc}
P_k&{\rm if}&j= k\\
0&{\rm if}&j\not= k
\end{array}\right.\\
2)\, {\rm completeness:}&\quad &\sum_{k=1}^n P_k={\bf 1}_n\,.
\end{eqnarray}

The non-singular linear operators of the set ${\rm Aut} ({\Bbb K}^n)\subset {\rm End} ({\Bbb K}^n)$  can be organized as a group with respect to  the matrix multiplication.  A specific feature here is the well-known rule $(AB)^{-1}=B^{-1} A^{-1}$ as in diagram (\ref{inverse}).
\begin{defin} 
The group ${\rm Aut} ({\Bbb K}^n)=GL(n, {\Bbb K})$  is called the group of general linear transformations of the vector space  ${\Bbb K}^n$. 
\end{defin}

\subsection{Matrix representations}

The above inspection of the matrix algebra suggests  all the abstract algebraic objects we consider could be {\em represented} by matrices. The advantage of this method is to reduce the study of the linear mappings to the simpler one of their matrices. 

\subsubsection{Representing vectors and tensors}

Any finite-dimensional vector space $V_n$ of dimension $n$ is isomorphic to the arithmetic space ${\Bbb K}^n$ such that  any vector $x\in V_n$ can be represented as a ket vector  of  same components of the arithmetic vector space ${\Bbb K}^n$  expanded  in the natural ket  basis of this space. Then the covectors of the vector space $V^*$ can be represented by the bra vectors of the space ${\Bbb K}^n$ having the same components.
\begin{defin}\label{defrep}
Given the pair of dual bases $(e,\hat e)$ in $(V_n,V_n^*)$ we say that the isomorphisms  ${\rm Rep}_e :V_n\to {\Bbb K}^n$  and ${\rm \hat Rep}_{\hat e} :V^*_n\to ({\Bbb K}^n)^*$ defined as 
\begin{equation}
{\rm Rep}_e (e_i)=|i\rangle\,,  \quad {\rm \hat Rep}_{\hat e} (\hat e^i)=\langle i|\,, \quad  i=1,2,...n,
\end{equation}\label{repres}
determine the (matrix) representation of the vector spaces $(V_n,V_n^*)$ corresponding to the mentioned dual bases. 
\end{defin}
Then  any vector $x=x^ie_i \in V_n$ is represented by the  vector  $|x\rangle={\rm Rep}_e (x)= x^i | i\rangle\in {\Bbb K}^n$  while any covector $\hat x=\hat x_i\hat e^i \in V^*_n$ is represented by the covector  $ \langle \hat x|={\rm \hat Rep}_{\hat e} (\hat x) \in ({\Bbb K}^n)^*$. Then the dual form reads $\langle \hat y,  x\rangle =\langle \hat y|x\rangle= {\rm \hat Rep}_{\hat e}(\hat y) {\rm Rep}_{e}(x)$.
\begin{theor}
Given a representation as in definition \ref{defrep} the dual isomorphisms $({\rm Rep}_e)^* : ({\Bbb K}^n)^*\to V^*_n$ and $({\rm \hat Rep}_{\hat e})^* : {\Bbb K}^n\to V_n$ satisfy
\begin{equation}
({\rm Rep}_e)^*=({\rm \hat Rep}_{\hat e})^{-1}\,, \quad ({\rm \hat Rep}_{\hat e})^*=({\rm Rep}_{e})^{-1}\,.
\end{equation}
\end{theor}
\begin{demo}
We write  $\langle \hat y , x\rangle=\langle \hat y | x\rangle ={\rm \hat Rep}_{\hat e}(\hat y){\rm Rep}_{e}(x)= \langle \hat y ,({\rm \hat Rep})^*_{\hat e}\circ {\rm Rep}_{e}(x)\rangle$ and deduce $({\rm \hat Rep})^*_{\hat e}\circ {\rm Rep}_{e}={\rm id}_V$ that prove the second equation since all the maps Rep are invertible as isomorphisms 
\end{demo}

The matrix representation of different classes of tensors is in general a complicated procedure which is  replaced often by the tensor calculus using indices. However there are some conjectures when  we need to represent tensors as matrices. Let us start with the set of vector spaces $V_i$ ($i=1,2,...N$), with ${\rm dim} V_i=n_i$, and the bases $e_{(i)}\subset V_i$ denoting with ${\rm Rep}_{e_{(i)}}  : V_i\to {\Bbb K}^{n_i}$ the representations corresponding to these bases.    
\begin{defin}
The tensor product of matrix representations is the representation 
\begin{equation}
{\rm Rep}_{e_{(1)}\otimes e_{(2)}...\otimes  e_{(N)}}={\rm Rep}_{e_{(1)}}\otimes {\rm Rep}_{e_{(2)}}....\otimes {\rm Rep}_{e_{(N)}}
\end{equation} 
corresponding to the bases $e_{(1)},\, e_{(2)},\,... e_{(N)}$.
\end{defin}
This representation is the isomorphism 
\begin{equation}
{\rm Rep}_{e_{(1)}\otimes e_{(2)}...\otimes  e_{(N)}}: V_1\otimes V_2...\otimes V_N \to
{\Bbb K}^{n_1}\otimes {\Bbb K}^{n_2}...\otimes {\Bbb K}^{n_N}
\end{equation} 
which maps each basis vectors $e_{(1)i_1}\otimes e_{(2)i_2}...\otimes  e_{(N)i_N}\in V_1\otimes V_2...\otimes V_N$ into the ket vectors of the tensor product basis in ${\Bbb K}^{n_1}\otimes {\Bbb K}^{n_2}...\otimes {\Bbb K}^{n_N}$ as
\begin{eqnarray}
&&{\rm Rep}_{e_{(1)}\otimes e_{(2)}...\otimes  e_{(N)}}(e_{(1)i_1}\otimes e_{(2)i_2}...\otimes  e_{(N)i_N})\\
&&~~~~~~~~~~={\rm Rep}_{e_{(1)}}(e_{(1)i_1})\otimes {\rm Rep}_{e_{(2)}}(e_{(2)i_2})....\otimes {\rm Rep}_{e_{(N)}}(e_{(N)i_N})\\
&&~~~~~~~~~~=|i_1\rangle \otimes |i_2\rangle...\otimes |i_N\rangle\,.  
\end{eqnarray} 
The ket vectors of the tensor product basis are denoted often by  $|i_1,i_2...i_N\rangle =|i_1\rangle \otimes |i_2\rangle...\otimes |i_N\rangle$. Obviously similar formulas hold for the dual spaces too.

This formalism can be used for  representing  the tensors of a given pair of vector spaces $(V,V^*)$ with ${\rm dim V}=n$. Thus the tensors of the space $T^k(V)$ can be represented in spaces $T^k({\Bbb K}^n)={\Bbb K}^n\otimes {\Bbb K}^n...=({\Bbb K}^n)^{\otimes k}$ while the dual tensors of the space $T_k(V)$ are represented in tensor products of dual spaces,   $T_k(({\Bbb K}^n))=({\Bbb K}^n)^*\otimes ({\Bbb K}^n)^*...=[({\Bbb K}^n)^*]^{\otimes k}$. 
Both these spaces can be seen as genuine arithmetic spaces whose bases may be constructed using the Kronecker product of matrices.
\begin{defin}
Given the $m_1\times n_1$ matrix $A$ and the $m_2\times n_2$ matrix $B$  their Kronecker product is the $m_1m_2\times n_1n_2$ matrix,
\begin{equation}\label{Kronprod}
A\otimes_{\kappa} B=\left(\begin{array}{cccc}
A^1_1B&A^1_2B&\cdots & A^1_{n_1}B\\
A^2_1B&A^2_2B&\cdots & A^2_{n_1}B\\
\vdots &\vdots &\cdots &\vdots \\
A^{m_1}_1B&A^{m_1}_2B&\cdots & A^{m_1}_{n_1}B\\
\end{array}\right)
\end{equation}
\end{defin}
Then next lemma considered for the simpler case of a tensor product of two arithmetic spaces offers us the key of the entire problem.
\begin{lem}
The tensor product of two arithmetic vector spaces, ${\Bbb K}^n\otimes {\Bbb K}^m$,  is isomorphic  with the space ${\Bbb K}^{n m}$ whose ket basis is the Kronecker product of the ket bases of the spaces  ${\Bbb K}^n$ and ${\Bbb K}^m$. 
\end{lem}
\begin{demo} The tensor product basis in ${\Bbb K}^n\otimes {\Bbb K}^m$ is formed by the vectors  $|i\rangle\otimes |j\rangle$. Then the natural isomorphism 
$\kappa : {\Bbb K}^n\otimes {\Bbb K}^m \to {\Bbb K}^{n m}$ so that  $\kappa(|i\rangle \otimes |j\rangle)=|(n-1)i+j\rangle\in {\Bbb K}^{n m}$ defines the Kronecker product $\otimes_{\kappa}=\kappa \circ \otimes$. Therefore, we can write  ${\Bbb K}^n\otimes_{\kappa} {\Bbb K}^m={\Bbb K}^{n m}$
\end{demo} Thus for any ket vectors $|x\rangle\in {\Bbb K}^n$ and $|y\rangle\in {\Bbb K}^m$, 
the vector $|x\rangle \otimes_{\kappa} |y\rangle$ is a ket vector too calculated according to the general rule given by Eq. (\ref{Kronprod}). 

Generalizing these results for tensor spaces we find the isomorphisms $T^k({\Bbb K}^n)\cong{\Bbb K}^{n^k}$ and $T_k({\Bbb K}^n)\cong ({\Bbb K}^*)^{n^k}$. 
Thus in the space $T^k({\Bbb K}^n)$ we can construct the ket basis formed by column matrices with $n^k$ components and corresponding bra bases of the same dimension in the dual vector space $T_k({\Bbb K}^n)$.  However, in the case of more complicated tensor spaces  the bases constructed using the Kronecker product are no longer ket or bra bases being formed by large matrices with many line and columns.  For example, the basis of the tensor space as $T^k_{\cdot\, l}({\Bbb K}^n)$ has to be formed by   $n^k\times n^l$ matrices. Obviously, such bases become useless in concrete calculations  such that one must consider only representations which do not mix  dual bases.  In fact,  the intuitive method of the Kronecker product remains of some  interest  only in very low dimensions but this  has applications in physics.

\subsubsection{Representing linear mappings}\label{Replinmap}

The principal application is the representation by matrices of  the linear mappings among different vector spaces defined over the same field ${\Bbb K}$. Let us start with two vector spaces $V$ and $W$ of dimensions $n={\rm dim}(V)$ and $m={\rm dim}(W)$ where we know the pairs of dual bases, $(e,\hat e)$ in $(V,V^*)$ and $(e',\hat e')$ in $(W,W^*)$. We consider the linear map   $f: V\to W$ and its dual map $f^*:W^*\to V^*$ given by definition (\ref{fstar}).  For any $x=x^i e_i \in V$ and $\hat y=\hat y_j \hat e'^j \in W^*$ we denote  $ y=f(x)=x^i f(e_i)\in W$ and $\hat x=f^*(\hat y)\in V^*$. Then we calculate the components in these bases defined by Eqs. (\ref{expv}) and (\ref{expcv}) that read
\begin{eqnarray}
&&y^i = \langle \hat e'^i,y\rangle =\langle \hat e'^i,f(x)\rangle=\langle \hat e'^i,f(e_j)\rangle x^j\,,\label{ff1}\\
&&\hat x_i = \langle \hat x,e_i\rangle=\langle f^*(\hat y), e_i\rangle =\hat y_k\langle f^*(\hat e'^k),e_i\rangle \,.\label{ff2}
\end{eqnarray} 
Now it is crucial to observe that the definition (\ref{fstar}) guarantees that  the above equations involve the same matrix elements
\begin{equation}\label{fij}
f^i_j =\langle f^*(\hat e'^i),e_j\rangle=\langle \hat e'^i,f(e_j)\rangle
\end{equation}
\begin{defin}\label{matf}
The $m\times n$ matrix $|f|$ with the matrix elements 
(\ref{fij})  is called the common matrix of  the maps $f$ and $f^*$  in the above considered bases. 
\end{defin} 
The problem is how to relate this matrix to a representation in arithmetic vector spaces. This can be done observing that  the commutative diagram 
\begin{equation}\label{dia1}
\begin{array}{rcl}
&f&\\
V~&\overrightarrow{\hspace*{20mm}}&W\\
\left.\begin{array}{r}
\\
{\rm Rep}_e\\
\\
\end{array} \right\downarrow~&&
\left\downarrow
\begin{array}{c}
\\
{\rm Rep}_{e'}\\
\\
\end{array}\right.\\
{\Bbb K}^n&\overrightarrow{\hspace*{20mm}}&{\Bbb K}^m\\
&{\rm Rep}_{(e^{\prime}, e)}(f)&
\end{array}
\end{equation}
allows us to associate the map $f$ to the linear map of the arithmetic spaces, 
\begin{equation}\label{repmat}
{\rm Rep}_{(e^{\prime}, e)}(f) ={\rm Rep}_{e^{\prime}}\circ f \circ{\rm Rep}_e^{-1} \in {\rm Lin}({\Bbb K}^n, {\Bbb K}^m)\,.
\end{equation}
\begin{theor}
In any representation the representation matrix coincides to  that of the represented map, i. e. ${\rm Rep}_{(e',e)}(f) =|f| $.  
\end{theor}
\begin{demo}
We calculate the matrix elements  
\begin{eqnarray}
\langle i'|{\rm Rep}_{(e',e)}(f) | j\rangle &=& \langle i'|{\rm Rep}_{e'}\circ f \circ{\rm Rep}_e^{-1}(|j\rangle)=\langle i'|({\rm Rep}_{e'}\circ f)(e_j)\nonumber\\
&=& \langle i'|{\rm Rep}_{e'}( f(e_j))=\langle \hat e'^i, f(e_j)\rangle\,,
\end{eqnarray} 
giving the desired result. 
\end{demo}
Thus we obtain the matrix representation of Eqs. (\ref{ff1}) and (\ref{ff2}),
\begin{equation}\label{ff12}
|y\rangle=|f|\,|x\rangle \,, \quad \langle \hat x|=\langle \hat y|\,|f|  \,,
\end{equation}
which is in accordance with our previous notations (\ref{AAstar}).

The representation theory helps us to study the properties of the linear mappings in the simple context of the theory of matrices and determinants. For example, many important results can be obtained starting with the next theorem.    
\begin{theor}
For any bases  $e\subset V_n$ and $e'\subset V_{m>n}$  there exists a linear map $f$ such that  $f(e_i)=e'_i$ for $i=1,2,...n<m$.
\end{theor} 
\begin{demo}
Since the basis $e$ is chosen we have the mapping ${\rm Rep}_e:V_n \to {\Bbb K}^n$. In addition we define the mapping $\phi:{\Bbb K}^n\to V_m$ so that $\phi(|x\rangle)=x^1 e'_1+x^2e'_2+...x^n e'_n$. Then we can put $f=\phi\circ {\rm Rep}_e$.  
\end{demo} 
Hereby one can prove that for any $f\in {\rm Lin}(V_n,V_m)$ there exists a pair of bases, $e\subset V_n$ and $e' \subset V_m$, such that the $m\times n$ representation matrix takes the canonic form
\begin{equation}
{\rm Rep}_{(e',e)}(f)=\left|
\begin{array}{cc}
{\bf 1}_{r\times r}&0\\
0&0
\end{array}\right|
\end{equation} 
where $r ={\rm dim}({\rm Im} f)={\rm rank} f$ is the {\em rank } of the map $f$. The principal consequence is that  ${\rm dim} ({\rm Im} f)+{\rm dim} ({\rm Ker} f)={\rm dim} (V_n)=n$.

Let us focus now on the endomorphisms of the vector space $V$ which form the algebra ${\rm End} (V)$
\begin{defin}
For each basis $e$ in $V$ the isomorphism  ${\rm Rep}_{(e,e)}: {\rm End} (V)\to {\rm End}({\Bbb K}^n)$ determines a matrix representation of the algebra ${\rm End} (V)$.
\end{defin}
In a given representations  the linear operators $f\in {\rm End}(V)$ are represented by the operators ${\rm Rep}_{(e,e)}(f)\in {\rm End}({\Bbb K}^n)$ given by Eq. (\ref{repmat}) but having now $n\times n$ matrices, ${\rm Rep}_{(e,e)}(f)=|f| $  with the matrix elements 
\begin{equation}\label{fffReee}
f^i_j=\langle i|{\rm Rep}_{(e,e)}(f)|j\rangle=\langle \hat e^i,f(e_j)\rangle=\langle f^*(\hat e^i),e_j\rangle\,.
\end{equation}
\begin{theor}\label{FG1} 
The following properties hold in any matrix representation of ${\rm End} (V)$.
\begin{eqnarray}
&&{\rm Rep}_{(e,e)}(f\circ g)={\rm Rep}_{(e,e)}(f)\,{\rm Rep}_{(e,e)}(g)\,,\label{prodmat} \\
&&{\rm  Rep}_{(e,e)}({\rm id}_V)={\bf 1}_{n}\, \quad n={\rm dim}(V)\,.
\end{eqnarray} 
\end{theor}
\begin{demo}
We denote $y^k=g^k_i x^i$ and $z^j=f^j_k y^k$ where $|f|$ and $|g|$ are the matrices representing the maps $f$ and $g$. Then we find $z^j=f^j_k g^k_ix^i$ which prove the first equation. The second one is obvious since ${\bf 1}_{n}$ is the matrix of the identity operator in ${\Bbb K}^n$ 
\end{demo}
Thus we understand that in a matrix representations the composition of maps reduces to simple  matrix multiplications. The natural consequence is:
\begin{cor}
For any invertible linear operator $f$ we have 
\begin{equation}\label{invmat}
{\rm Rep}_{(e,e)}(f^{-1})={\rm Rep}_{(e,e)}(f)^{-1}
\end{equation}
\end{cor}
The dual map $f^*:V^*\to V^*$  of the  map of $f: V\to V$ has the same matrix $|f|$ having the matrix elements $\langle f^*(\hat e^i), e_j\rangle=\langle \hat e^i, f(e_j)\rangle$. This has the direct and dual actions as given by Eq. (\ref{ff12}). 

Another advantage of the representation theory  is that we can apply the theory of determinants.
\begin{defin}
The ${\Bbb K}$-number ${\rm det}(f)={\rm det}|{\rm Rep}_{(e,e)}(f)|$, is called the determinant of the linear operator $f$.
\end{defin} 
The linear operators $f\in {\rm End}(V)$ with ${\rm det}(f)=0$ are degenerate operators (or projectors) having ${\rm Ker} f \not= \{0\}$.  The  non-degenerate  linear operators are invertible and form the group ${\rm Aut}(V)$.
\begin{defin}
The map ${\rm Rep}_{(e,e)}: {\rm Aut} (V)\to {\rm Aut}({\Bbb K}^n)= GL(n,{\Bbb K})$ is a
matrix representation of the group ${\rm Aut}(V)$.
\end{defin}

Finally, we note that the method of matrix representation can be applied in the same manner to the multilinear maps.  For example, starting with two operators $f$ and $g$ represented by their matrices $|f|\in {\rm End}({\Bbb K}^n)$ and $|g|\in {\rm End}({\Bbb K}^m)$ we can represent the map $f\otimes g$ either as a tensor product $|f|\otimes |g| \in {\rm End}({\Bbb K}^n \otimes{\Bbb K}^m)$ or as a Kronecker product $|f|\otimes_{\kappa} |g| \in {\rm End}({\Bbb K}^{nm})$. Both these methods are used for solving concrete problems in physics.

\subsubsection{Changing the representation}

The operators  $f\in {\rm End}(V)$ can be represented by matrices in any basis of $V$. Let us see now what happens when the bases of the vector spaces $(V,V^*)$ are changed $(e,\hat e)\to (e',\hat e')$ changing thus the representation. The change of basis means to change simultaneously the basis and the components such that the vectors remain the same, i. e. $x=x^ie_i=x^{\prime\, i}e'_i$. Then the ket vector $|x\rangle \in {\Bbb K}^n$ representing $x$ is transform into the new ket vector $|x'\rangle $ taking into account that the natural basis of   ${\Bbb K}^n$ can not be changed.
This conjecture corresponds to the commuting  diagram 
\begin{equation}\label{dia2}
\begin{array}{rcl}
&{\rm id}_V&\\
V~&\overrightarrow{\hspace*{20mm}}&V\\
\left.\begin{array}{c}
\\
{\rm Rep}_e\\
\\
\end{array} \right\downarrow~&&
\left\downarrow
\begin{array}{c}
\\
{\rm Rep}_{e'}\\
\\
\end{array}\right.\\
{\Bbb K}^n&\overrightarrow{\hspace*{20mm}}&{\Bbb K}^n\\
&M&
\end{array}
\end{equation}
Then we denote by $x^i$ and $\hat y_i,$ the components of a vector $x$ and a covector $\hat y$ in the bases $(e,\hat e)$ and by $x'^i$ and $\hat y'_i$ the similar components in the new bases $(e',\hat e')$.  
\begin{theor}
Changing the bases the components of the vectors and covectors transform as
\begin{equation}
|x'\rangle=M\,|x\rangle \Rightarrow x^{\prime \, i}=M^i_j x^j\,,\quad \langle \hat y'|=\langle \hat y|\,M^{-1}\Rightarrow \hat y'_i=(M^{-1})^j_i\hat y_j\,,
\end{equation} 
where  $M\in {\rm End}({\Bbb K}^n)$ is the matrix
\begin{equation}\label{repmat1} 
M={\rm Rep}_{(e',e)}({\rm id}_V)= {\rm Rep}_{e'} \circ{\rm Rep}_e^{-1} \in {\rm Aut}({\Bbb K}^n)
\end{equation}
called  the change of basis operator. 
\end{theor}
\begin{demo} 
The operator $M$ and its dual operator, $M^*={\rm \hat Rep}_{\hat e}\circ {\rm \hat Rep}^{-1}_{\hat e'}\in {\rm Aut}[({\Bbb K}^n)^*]$, act as in Eq. (\ref{AAstar}).  Then we calculate $M\, |x\rangle =(M\circ {\rm Rep}_e)(x)={\rm Rep}_{e'}(x)=|x'\rangle$ and 
similarly $\langle \hat y'|\,M=M^*(\langle \hat y '|)={\rm \hat Rep}_{\hat e}(\langle \hat y|)=\langle \hat y|$ giving $\langle \hat y'|=\langle \hat y|\,M^{-1}$
\end{demo}  
Conversely,  any non-singular linear operator $M\in {\rm Aut}({\Bbb K}^n)$ transforms the basis $e$ into the new one, $e'$, such that ${\rm Rep}_{e'}=M\circ {\rm Rep}_e$. Then we find the new dual  bases in $(V,V^*)$ given by this transformation, 
\begin{equation}\label{eeM} 
e'_{i}=e_j (M^{-1})^j_i \,, \quad \hat e'^{i}=\hat e^j M^i_j\,. 
\end{equation}
As was expected the dual form remains invariant under the change of representation, $\langle \hat y'|x'\rangle =\langle \hat y|x \rangle$. 
\begin{cor}
Given two arbitrary pairs of dual bases  $(e,\hat e)$ and $(e',\hat e')$ these are related by Eqs. (\ref{eeM}) where
\begin{equation}
M^i_j=\langle \hat e'^i, e_j\rangle\,,  \quad (M^{-1})^j_i=\langle \hat e^j, e'_i\rangle\,.
\end{equation}
\end{cor}
\begin{demo}
We calculate $M^i_j=\langle i'| M|j\rangle=\langle i'|{\rm  Rep}_{e'}\circ {\rm Rep}_e^{-1}(|j\rangle)=\langle i'|{\rm Rep}_{e'}(e_j) =\langle \hat e'^i,e_j\rangle$ and similarly 
$(M^{-1})^j_i=\langle j|M^{-1}|i'\rangle= \langle \hat e^j, e'_i\rangle$ 
\end{demo}
Under a basis transformation the tensor components transform as
\begin{equation}
\tau^{i_1 i_2...}_{\,. \,\, .\,\dots  j_1 j_2...}\to \tau^{\prime\, \,i_1 i_2...}_{\,\,\,. \,\, .\,\dots  j_1 j_2...}=M^{i_1}_{k_1} M^{i_2}_{k_2}(M^{-1})^{l_1}_{j_1} (M^{-1})^{l_2}_{j_2}\cdots\tau^{k_1 k_2...}_{\,. \,\, .\,\dots  l_1 l_2...}  \,.
\end{equation}

Finally, let us see how the operator matrices transform when we change the representation.
\begin{theor}
Given $f\in {\rm Aut}(V)$ and two different  bases $e$ and $e'$ in $V$  then the representation matrices in these bases satisfy 
\begin{equation}\label{FMFM} 
{\rm Rep}_{(e',e')}(f)={\rm Rep}_{(e',e)}({\rm id}_V)\,{\rm Rep}_{(e,e)}(f)\,{\rm Rep}_{(e,e')}({\rm id}_V)
\end{equation}
where ${\rm Rep}_{(e,e')}({\rm id}_V)={\rm Rep}_{(e',e)}({\rm id}_V)^{-1}$.
\end{theor}
\begin{demo}
The last equation can be obtained  straightforwardly from Eqs.  (\ref{repmat1}) and (\ref{invmat}). Using then  Eq. (\ref{repmat})  we calculate 
\begin{equation}
{\rm Rep}_{(e',e)}({\rm id}_V)\circ {\rm Rep}_{(e,e)}(f)\circ {\rm Rep}_{(e,e')}({\rm id}_V)={\rm Rep}_{e'}\circ  f\circ  {\rm Rep}_{e'}^{-1}
\end{equation}
which leads to the desired result using Eq. (\ref{prodmat}).
\end{demo} In the usual (self-explanatory) notation Eq. (\ref{FMFM}) reads 
\begin{equation}
|f'|=M\,| f|\,M^{-1}\quad \Rightarrow \quad 
f^{\prime\,i}_j=M^i_k\,f^k_l(M^{-1})^l_j\,.
\end{equation}
This equation shows that ${\rm det} f'= {\rm det} f$ which proves  the following theorem.  
\begin{theor}
The  determinant of the linear operators are invariants, independent on the choice of the basis in which these are calculated. 
\end{theor} 
Then the theorem (\ref{detAB}) gives us the rule  ${\rm det}(f\circ g)={\rm det}(f)\,{\rm det}(g)$ that holds in any representation.

\section{Inner product vector spaces}

The natural duality is not enough for giving rise to invariant quantities with some physical meaning. These are  produced in real or complex vector spaces equipped with definite or indefinite inner products.

\subsection{Unitary vector spaces}

The unitary (or inner product) spaces are largely used in quantum mechanics as spaces carrying various unitary finite-dimensional representations (of group or algebras) that may be of physical interest. In what follows we focus on such spaces and their associated operators.

\subsubsection{Topologies in vector spaces}

\begin{defin}
A topological vector space is a vector space with a Hausdorff topology in which the addition and scalar multiplication are  continuous.
\end{defin}
A convenient topology of open sets can be introduced with the help of real quantities associated with each vector $x\in V$ and compatible with the operations of $V$. Let us consider complex vector spaces $V,\, W,...$ defined over the field ${\Bbb C}$.
\begin{defin}\label{defnorm}
Suppose $V$ is a complex vector space.  We say that the map $\|~\|:V\to {\Bbb R}$ is a norm if:

1)  $\| x\| \ge 0$ and $ \|x\|= 0$ if and only if $x=0\in V$,

2)  $\|x+y\|\le \{x\|+\|y\|$,

3)  $\|\alpha x\|=|\alpha|\|x\|$ for all $\alpha \in {\Bbb C}$.
\end{defin}
The norm can give a distance $d(x,y)=\|x-y\|$ between any vectors $x,\,y\in V$ tanks to the fact that $V$ is a vector space.   
\begin{defin}
We say that $V$ is a normed (vector) space if this has the topology generated by the family of open balls $B(x,r)=\{y\in V| \,\|x-y\|<r\}$. 
\end{defin}
The next step is to define a quantity able to generate the norm as in the familiar case of the three-dimensional real-valued vectors. This is the inner product that can be defined in complex or real vector spaces.

When we work with complex vector spaces we say that a map $f:V\to W$ which satisfy $f(x+y)=f(x)+f(y)$ is linear if $f(\alpha x)=\alpha f(x)$ but {\em anti-linear} whether $f(\alpha x)=\overline{\alpha} f(x)$.
\begin{defin}
The map $H: V\times V\to {\Bbb C}$ is a Hermitian form (H-form) on the complex  vector space $V$  if
\begin{eqnarray}
&1)&H(x,\alpha y +\beta y')=\alpha\, H(x,y)+\beta\, H(x,y')\,, \quad \forall\, x,y,y' \in V\,, ~ \alpha,\beta \in {\Bbb C}\nonumber\\
&2)& H(y,x)=\overline{H(x,y)}\nonumber
\end{eqnarray}
\end{defin}
This is linear in the second argument but anti-linear in the first one since 
$H(\alpha x +\beta x',y)=\overline{\alpha}\,H(x,y)+\overline{\beta}\,H(x',y)$ as it results combining the above axioms. Moreover, from the last one it results that for $x=y$ the H-form  is a real number, $H(x,x)\in {\Bbb R}$.
\begin{defin}\label{H-form}
The H-form $h$, is {non-degenerate} if its (left) kernel (which is equal to the right one) reduces to the set $\{0\}$,
\begin{equation}
{\rm Ker}\, h =\{x\, | H(x,y)=0\,, \forall \, y\in V\}= \{y\, | H(x,y)=0\,, \forall \, x\in V\}=\{0\}\,.
\end{equation} 
\end{defin}
\begin{defin}
A non-degenerate H-form $(~,~):V\times V\to {\Bbb C}$ is called a inner product if satisfies 
the axioms of definition (\ref{H-form}),  

1) $ (x,\alpha y +\beta y')=\alpha\, (x,y)+\beta\, (x,y')$ for all $x,y,y' \in V\,, ~ \alpha,\beta \in {\Bbb C}$

2) $(y,x)=\overline{(x,y)}$\\
and complies with the specific axiom

{\rm 3)} $(x,x)\ge 0$ and $(x,x)=0$ if and only if $x=0\in V$.
\end{defin}
This last axiom guarantees that the inner product can define the associate {\em norm}  $\|x\|=\sqrt{(x,x)}$ that satisfies the conditions of definition (\ref{defnorm}). In addition  there is a remarkable supplemental property.
\begin{theor}
The norm derived from a inner product satisfy the Cauchy-Bunyakovsky-Schwarz inequality, 
$|(x,y)|\le \|x\|\,\|y\|$.
\end{theor}
\begin{demo}
Let $x=(x,y)(y,y)^{-1}y+z $ and verify that $z$ is orthogonal to $y$, i. e. $(z,y)=0$. Then we find \begin{equation}
\|x\|^2=(x,x)=\left | \frac{(x,y)}{(y,y)}\right | \|y\|^2+\|z\|^2=\frac{|(x,y)|^2}{\|y\|^2}+\|z\|^2\ge \frac{|(x,y)|^2}{\|y\|^2}
\end{equation}
which proves the inequality \end{demo}
\begin{defin}
A finite-dimensional complex vector space  equipped with a inner product is called a unitary vector space or inner product space. 
\end{defin}
The inner product determines orthogonality properties in $V$.  We recall that two vectors $x,y\in V$ are orthogonal if $(x,y)=0$ while  two subspaces $V_1$ and $V_2$ are orthogonal if  $(x_1,x_2)=0$ for all $x_1\in V_1$ and $x_2\in V_2$. This is denoted as $V_1\perp V_2$. Each subspace $V_1\subset V$ has an orthogonal complement $V_2\perp V_1$ so that  $V=V_1\oplus V_2$. This means that any $x\in V$ can be projected on these subspaces as $x=x_1+x_2$ such that $x_1\in V_1$ and $x_2\in V_2$ are uniquely determined.

\subsubsection{Orthonormal dual bases}

Now we have two types of forms, the dual form $\langle~,~\rangle:V^*\times V \to {\Bbb C}$ and the inner product $(~,~):V\times V\to {\Bbb C}$. These can be naturally related by an anti-linear mapping $\phi_R$ defined by the Riesz representation theorem.
\begin{theor}{\rm (Riesz)}\label{Riesz}
Let $V$ be a unitary vector space and $V^*$  its dual.  There exists a bijective anti-linear map $\phi_R:V\to V^*$ such that    
\begin{equation}\label{psiR}
\langle \phi_R(x),y\rangle=(x,y)\,.
\end{equation}
\end{theor}
\begin{demo}
Let us take  $x\in {\rm Ker}\, \phi_R=\{x|\phi_R(x)=0\}$ which implies $(x,x)=0$ but then $x=0\in V$. Therefore ${\rm Ker}\, \phi_R=\{0\}$ and $\phi_R$ is bijective.  Furthermore, we derive
$\langle \phi_R(\alpha x),y\rangle=(\alpha x,y)=\overline{\alpha} (x,y)=\langle \overline{\alpha}\phi_R(x), y\rangle$  which means that $\phi_R(\alpha x)=\overline{\alpha}\phi_R(x)$
\end{demo} The principal consequence is that $V^*$ is a unitary space too.
\begin{cor}
Suppose $V$ is a unitary vector space with the inner product $(~,~)$. Then its dual space $V^*$ is a unitary space with respect to the inner product $(~,~)_*: V^*\times V^*\to {\Bbb C}$ defined as
\begin{equation}
(\hat x, \hat y)_*=\langle \hat y,\phi_R^{-1}(\hat x)\rangle=(\phi_R^{-1}(\hat y),\phi_R^{-1}(\hat x)) \,,\quad \hat x,\,\hat y \in V^*\,.
\end{equation}   
\end{cor}
\begin{demo}
This definition gives a linear form in $\hat y$ and anti-linear in $\hat x$
\end{demo} 
\begin{rem}\label{phiself}
According to definition (\ref{selfdual}) the Riesz isomorphism is selfdual, $\phi_R^*=\phi_R$. 
\end{rem}

In the unitary spaces $(V,V^*)$ we consider the pair of dual bases $(e,\hat e)$.
\begin{defin}
The dual bases $(e,\hat e)$ are orthonormal if  the Riesz  isomorphism $\phi_R$ acts as  
\begin{equation}
\phi_R(e_i)=\hat e^i
\end{equation}
for all $e_i\in e$.  
\end{defin}
\begin{lem}
The vectors of a pair of orthonormal dual bases satisfy
\begin{equation}
(e_i, e_j)=\delta_{ij}\,, \quad (\hat e^i,\hat e^j)_*=\delta_{ij}\,.
\end{equation} 
\end{lem}
\begin{demo}
$(e_i,e_j)=\langle \phi_R(e_i),e_j\rangle=\langle \hat e^i, e_j\rangle=\delta_{ij}$ and $(\hat e^i,\hat e^j)_*=(\phi_R^{-1}(\hat e^j),\phi_R^{-1}(\hat e^i))=(e_j,e_i)=\delta_{ij}$. 
\end{demo} Note that, in general,  the Riesz isomorphism gives 
\begin{equation}
\phi_R (x)=\phi_R(x^i e_i)=\overline{x^i}\phi_R( e_i)=\sum _i\overline{x^i}\,\hat e^i\,.
\end{equation}  
\begin{rem} The isomorphism $\phi_R$ could hide or dissimulate  the real position of indices  in some formulas  forcing us to give up the Einstein convention using explicit sums in order to avoid confusion.
\end{rem}
Thus in  the orthonormal dual bases $(e, \hat e)$ we have
\begin{equation}
x=x^i e_i =\sum_i e_i (e_i,x)\,,
\end{equation} 
since now $x^i=\langle \hat e^i, x\rangle=\langle\phi_R(e_i), x\rangle=(e_i,x)$.
The inner product can also be expanded in this basis, 
\begin{equation}
(x,y)=(x^i e_i,y^j e_j)=\overline{x^i} y^j (e_i,e_j)=\sum_i \overline{x^i} y^i=\sum_i (x,e_i)(e_i,y)\,,
\end{equation} 
using explicit sums. In particular, the norm reads $\|x\|= (\sum_i |x^i|^2)^{\frac{1}{2}}$.

The orthonormal bases can be changed according to the rule given by Eq. ({\ref{eeM}}). In general, an arbitrary matrix $M$ can change orthonormal bases in other bases that may be not orthonormal.
\begin{theor}\label{ortbasM}
If  $(e,\hat e)$ and $(e',\hat e')$ are orthonormal bases then their transformation matrix $M$ is unitary, $M^+=M^{-1}$.  
\end{theor}
\begin{demo}
We take $\delta_{i,k}=(e'_i,e'_k)=\overline{(M^{-1})^j_i}(M^{-1})^l_k\,(e_j,e_l)=\sum_j \overline{(M^{-1})^j_i}(M^{-1})^j_k$ that implies $M^+=M^{-1}$
\end{demo}
\begin{cor}
The Riesz isomorphism is independent on the choice of the orthonormal bases.
\end{cor}
\begin{demo}
We consider the bases of the above theorem and write $\phi_R (e'_i)= \overline{(M^{-1})^j_i}\phi_R(e_j)=M^i_j \phi_R(e_j)$ since $M$ is unitary. Then  
$\phi_R(e_j)=\hat e^j \, \Rightarrow \phi_R (e'_i)=M^i_j \hat e^j=\hat e'^i$
\end{demo}

\subsubsection{Operators on unitary vector spaces}

The operators $f\in {\rm End} (V)$  defined on a unitary vector space have supplemental properties related to the inner product. 
\begin{defin}
Given the dual operators $f:V\to V$ and  $f^*:V^*\to V^*$ we say that: 

1) the operator $f^{\dagger}:V\to V$ is the adjoint operator of $f$  if  the identity  $(f^{\dagger}(x), y)=(x, f(y))$ holds for any $x,y \in V$,

2) the operator $(f^*)^{\dagger}:V^*\to V^*$ is the adjoint operator of  $f^*$  if  we have  $((f^*)^{\dagger}(\hat x), \hat y)_*=(\hat x, f^*(\hat y))_*$  for any $\hat x,\hat y \in V^*$ .
\end{defin}
The following calculation rules are obvious: $(f^{\dagger})^{\dagger}=f$ and
\begin{equation}
(f+g)^{\dagger}=f^{\dagger}+g^{\dagger}\,,\quad (f\circ g)^{\dagger} =g^{\dagger}\circ f^{\dagger}\,, \quad  (\alpha f)^{\dagger}=\overline{\alpha}f^{\dagger}\,,\,\,\alpha\in {\Bbb C}\,,
\end{equation} 
and similarly for the dual operators.
\begin{theor}
The dual operators $f$ and $f^*$ and their adjoint operators $f^{\dagger}$ and  $(f^*)^{\dagger}$ are related through
\begin{equation}\label{ffstarcr}
f= \phi_R^{-1}\circ (f^*)^{\dagger}\circ \phi_R\,, \quad f^{\dagger}= \phi_R^{-1}\circ f^*\circ \phi_R \, \Rightarrow (f^{\dagger})^*=(f^*)^{\dagger}\,.
\end{equation}
\end{theor}
\begin{demo}
$\langle \hat y, f\circ\phi^{-1}_R(\hat x)\rangle=(\hat x,f^*(\hat y))_*=((f^*)^{\dagger}(\hat x), \hat y)_*=\langle \hat y, \phi_R^{-1}\circ(f^*)^{\dagger}(\hat x)\rangle$ provs the first formula. For the second one we thake $(f^{\dagger}(x),y)=(x,f(y))=\langle\phi_R(x),f(y)\rangle=\langle f^*\circ\phi_R(x),y\rangle=(\phi_R^{-1}\circ f^*\circ\phi_R(x),y)$. Then we change $f\to f^{\dagger}$ in the last formula
\end{demo} Thus we show  that the operations dual  $ ^*$ and adjoint $ ^{\dagger}$ {\em commute}. These results correspond to the commutative diagrams:
\begin{equation}\label{dia4}
\begin{array}{rcl}
&f&\\
V~&\overrightarrow{\hspace*{20mm}}&V\\
\left.\begin{array}{c}
\\
\phi_R\\
\\
\end{array} \right\downarrow~&&
\left\downarrow
\begin{array}{c}
\\
\phi_R\\
\\
\end{array}\right.\\
V^*&\overrightarrow{\hspace*{20mm}}&V^*\\
&(f^*)^{\dagger}&
\end{array}
\begin{array}{rcl}
&f^{\dagger}&\\
V~&\overrightarrow{\hspace*{20mm}}&V\\
\left.\begin{array}{c}
\\
\phi_R\\
\\
\end{array} \right\downarrow~&&
\left\downarrow
\begin{array}{c}
\\
\phi_R\\
\\
\end{array}\right.\\
V^*&\overrightarrow{\hspace*{20mm}}&V^*\\
&f^*&
\end{array}
\end{equation}
This conjecture suggests how the adjoint operators of the operators of ${\rm Lin}(V,V^*)$ could be defined.
\begin{defin}
The adjoint operator of the operator $f:V\to V^*$ is the operator $f^{\dagger}:V^*\to V$ which satisfies
\begin{equation}
(x, f^{\dagger}\circ\phi_R(y))=(\phi_R^{-1}\circ f(x),y)\,.
\end{equation}
\end{defin}
\begin{rem}\label{ciuciubac}
According to this definition we can write $\phi_R^{\dagger}=\phi_R^{-1}$.
\end{rem}
Notice that $\phi_R$ with this property may not be interpreted as a unitary operator since these must be exclusively operators of ${\rm End}(V)$ or  ${\rm End}(V^*)$. However, this property is a good tool in current calculations. For example, Eqs.  (\ref{ffstarcr}) can be seen now as being adjoint to each other. 

The matrices of these operators calculated  in orthonormal dual bases  according to the general rules discussed in section (\ref{Replinmap}) are
\begin{eqnarray}
f^i_j&=&\langle \hat e^i, f(e_j)\rangle =\langle \phi_R(e_i),e_j\rangle =(e_i,f(e_j))\,, \label{cucu1}\\
(f^{\dagger})^i_j&=&(e_i,f^{\dagger}(e_j))=\overline{(f^{\dagger}(e_j),e_i)}=\overline{(e_j,f(e_i))}=\overline{f_{i}^j}\,.\label{cucu2}
\end{eqnarray}
Thus we recover the well-known result that in orthonormal bases the matrices of the adjoint operators are related through the Hermitian conjugation, $|f^{\dagger}|=|f|^+$. On the other hand, we see that the map $\phi_R$ hides the index position in the last term of Eq. (\ref{cucu1}). For this reason we expand functional equations of the form
\begin{equation}
y=f(x) \quad \Rightarrow \quad y^i= f_{j}^i x^j=\sum_j (e_i,f(e_j))(e_j,x)
\end{equation}
using explicit sums.
\begin{defin}
The equation $f(x_{\lambda})=\lambda x_{\lambda}$ is called the eigenvalue problem of the operator $f\in {\rm End}(V)$. The solutions are the eigenvectors $x_{\lambda}$ corresponding to the eigenvalues $\lambda$.
\end{defin}
The eigenvector $x=0\in V$ is always a solution of this equation for any $\lambda$. We say that this is a {\em trivial} eigenvector. For $\lambda=0$ we can meet two cases: 1) there is only the trivial eigenvector $x=0$ when the operator is non-singular and ${\rm Ker} f =\{0\}$, and 2) there are non-trivial eigenvectors that form the non-trivial subspace $V_0={\rm Ker} f $ of the singular operator $f$. In the first case we say that $\lambda=0$ is a {\em trivial} eigenvalue. 
\begin{lem}
The set of the eigenvectors corresponding to a non-trivial eigenvalue $\lambda$ form a non-trivial subspace, $V_{\lambda}\subset V$, called the eigenspace of $\lambda$. 
\end{lem}
\begin{demo}
If $\lambda=0$ is non-trivial then we are in the case 2) when the corresponding subspace $V_0$ is non-trivial.  For $\lambda\not=0$ we see that any linear combinations of eigenvectors corresponding to $\lambda$ is an eigenvector too 
\end{demo}
\begin{lem}
The set of non-trivial eigenvalues,  ${\rm Sp}(f)$, called  the spectrum of $f$,  is formed by the solutions of the secular equation, 
\begin{equation}\label{secular}
{\rm det}\left(|f|-\lambda {\bf 1}_n\right)=0\,, \quad n={\rm dim} V\,.
\end{equation}
\end{lem}
\begin{demo}
The eigenvalue problem can be expanded in orthonormal bases obtaining the homogeneous system
 \begin{equation}\label{homog}
\left(f_{k}^i-\lambda \delta_{k}^i\right)x_{\lambda}^k=0\,,
\end{equation}
which has non-trivial solutions only when   its determinant vanishes
\end{demo}
\begin{defin}
We say that an operator which satisfies $f^{\dagger}=f$ is a selfadjoint operator.
\end{defin}
\begin{lem}
The dual operator $f^*$ of a selfadjoint operator $f$ is a selfadjoint operator, $(f^*)^+=f^*$, which satisfies $f^*=\phi_R f \phi_R^{-1}$.  
\end{lem}
\begin{demo} We replace $f^+=f$ in Eq. (\ref{ffstarcr}b) \end{demo} In orthonormal bases the matrices of the selfadjoint operators are Hermitian, $|f^*|=|f|=|f|^+$ as it results from Eqs. (\ref{cucu1}) and (\ref{cucu2}).
\begin{theor}{\rm (Spectral theorem)}\label{spectheor}
A selfadjoint operator has only real eigenvalues and the eigenvector subspaces corresponding to different non-trivial eigenvalues are orthogonal.
\end{theor}
\begin{demo}
Let us consider two eigenvectors of the operator $f$ that obey $f(x_{\lambda_1})=\lambda_1 x_{\lambda_1}$ and $f(x_{\lambda_2})=\lambda_2 x_{\lambda_2}$ and calculate
$(x_{\lambda_1},f(x_{\lambda_2})) =\lambda_2 (x_{\lambda_1},x_{\lambda_2})$ and 
$(f(x_{\lambda_1}),x_{\lambda_2}) =\overline{\lambda_1} (x_{\lambda_1},x_{\lambda_2})$.
But $f=f^{\dagger}$ which means that $(x_{\lambda_1},f(x_{\lambda_2}))=(f(x_{\lambda_1}),x_{\lambda_2})$ and, therefore, we obtain the identities
\begin{equation}
(\lambda_2-\overline{\lambda_1})(x_{\lambda_1},x_{\lambda_2})=0\,,\quad \forall \lambda_1,\lambda_2\in {\rm Sp}(f)\,.
\end{equation}
Taking first $\lambda_1=\lambda_2$ and $x_{\lambda_1}=x_{\lambda_2}$ we find 
$(\lambda_1-\overline{\lambda_1})(x_{\lambda_1},x_{\lambda_1})=0$ that implies 
$\lambda_1=\overline{\lambda_1}$ since $(x_{\lambda_1},x_{\lambda_1})\not=0$ as long as $x_{\lambda_1}$ is a non-trivial eigenvector. Otherwise, if $\lambda_1\not=\lambda_2$, we must have $(x_{\lambda_1},x_{\lambda_2})=0$ such that $V_{\lambda_1}\perp V_{\lambda_2}$
\end{demo}

This result is important for analysing the action of the selfadjoint operators  $f=f^{\dagger}\in {\rm End} (V)$ using algebraic methods. First of all one has to derive the spectrum ${\rm Sp}(f)=\{\lambda_1,\lambda_2,...\lambda_k\}$ by solving the secular equation (\ref{secular}). Furthermore, one must identify the eigenspaces $V_{\lambda_i}$ solving the homogeneous system (\ref{homog}) for each of the non-trivial eigenvalues.  One obtain thus the orthogonal decomposition
\begin{equation}\label{VVV}
V=V_{\lambda_1}\oplus V_{\lambda_2}...\oplus V_{\lambda_k}
\end{equation}   
where the subspaces $V_{\lambda_i}$ have dimensions $n_i$ such that $n_1+n_2...+n_k=n$. In each eigenspace $V_{\lambda_i}$  one may choose an orthonormal basis, 
\begin{equation}
\epsilon_{(i)}=\{\epsilon_{(i)\,j}|(\epsilon_{(i)\,j},\epsilon_{(i)\,l})=\delta_{jl}\,, j,l=1,2...n_i\}\subset V_{\lambda_i} 
\end{equation}
and construct  the whole orthonormal basis $\epsilon =\epsilon_{(1)}\cup \epsilon_{(2)}...\cup \epsilon_{(k)}$ of the entire space $V$. 
\begin{cor}
For any selfadjoint operator $f$ there exists an orthonormal basis $\epsilon$ in which its matrix take the canonical form 
\begin{equation}\label{canform}
|f|_{\epsilon} ={\rm diag}(\underbrace{\lambda_1,\lambda_1...,\lambda_1}_{n_1}, \underbrace{\lambda_2,\lambda_2...,\lambda_2}_{n_2},... \underbrace{\lambda_k,\lambda_k...,\lambda_k}_{n_k})  
\end{equation}
\end{cor}
\begin{demo}
This is just the basis $\epsilon$ introduced above whose vectors are orthogonal among themselves and satisfy $f(\epsilon_{(i)\, j})=\lambda_i \epsilon_{(i)\,j}$. Then we have   
$(\epsilon_{(i)\, l},f(\epsilon_{(i)\, j}))=\lambda_i \delta_{lj}$ for $i=1,2...k$
\end{demo}
\begin{defin}
The numbers $n_i={\rm dim}(V_{\lambda_i})$ are the multiplicities of the eigenvalues $\lambda_i$.
\end{defin}
Note that one of eigenvalues, say $\lambda_1$, can be zero. Then its multiplicity $n_1$ gives  
${\rm dim}({\rm Ker} f)=n_1>0$ which shows that $f$ is singular.
\begin{defin}
The selfadjoint operator $p \in {\rm End}(V)$ that satisfies $p\circ p=p$ is  called the projection operator. The subspace $p(V)={\rm Im}\, p\in V$ is the projection subspaces.
\end{defin}
The trivial case is when $p={\rm id}_V$. In general, $p$ is singular so that there exists an orthonormal basis in which his matrix takes the diagonal form,
\begin{equation}
|p|={\rm diag}(\underbrace{1,1,...,1}_{n_p},\underbrace{0,0,...,0}_{n-n_p})\,,
\end{equation}
allowed by the condition $|p|^2=|p|$. The projection space ${\rm Im}\, p$ has the dimension $n_p$ while  the dimension of ${\rm Ker}\, p$ is $n-n_p$. Since $p$ is selfadjoint these subspace are orthogonal so that $V={\rm Im}\, p \oplus {\rm Ker}\, p$. The complementary projection operator of $p$ is the operator ${\rm id}_V -p$ for which ${\rm Im}({\rm id}_V -p) ={\rm Ker}\, p$.
\begin{defin}\label{ortproj}
A system of projection operators $p_i,i=1,2...n$ is called orthogonal if 
\begin{equation}
p_i\circ p_j =\left\{\begin{array}{ccc}
p_i&{\rm if}& i=j\\
0 &{\rm if} & i\not=j
\end{array}\right.
\end{equation}
If, in addition, we have $p_1+p_2...+p_n={\rm id}_V$ then we say that the system of projection operators is complete.  
\end{defin}
The simplest example is of the above orthogonal projection operators $p$ and ${\rm id}_V-p$ that form a complete system.  A more general example we met in section (\ref{Aritspaces}).

Important applications are the spectral representations of the selfadjoint operators $f$ for which we know the spectrum ${\rm Sp}(f)$ and the corresponding eigenspaces, $V_{\lambda}$. Then we define the projection operators on these eigenspaces, $p_{\lambda}$, which give the projections $p_{\lambda}(V)=V_{\lambda}$.
\begin{lem}
The projectors $\{p_{\lambda}|\lambda\in {\rm Sp}(V)\}$ form a complete system of orthogonal projectors  that give the spectral representation
\begin{equation}\label{spreprep}
f=\sum_{\lambda\in {\rm Sp}(f)} \lambda\, p_{\lambda}
\end{equation}  
\end{lem}
\begin{demo}
The projection operators $p_{\lambda}$ and $p_{\lambda'}$ are orthogonal when   $\lambda'\not=\lambda$ since then $V_{\lambda} \perp V_{\lambda'}$. The system is complete corresponding to the orthogonal sum (\ref{VVV}). In the orthonormal basis $\epsilon$ this representation gives the canonic form (\ref{canform}) 
\end{demo} The advantage of the spectral representations is that these are basis-independent.

\subsubsection{Isometries}

The vector spaces can be related through linear transformations that preserve the algebraic properties. Now we consider unitary spaces each one having its own inner product. Since we shall handle many spaces in this section we denote by $(~,~)_V$ the inner product of a unitary space $V$.
\begin{defin}\label{H-form}
The linear map $f:V\to W$ between two unitary vector spaces is an isometry if this
preserves the inner product, $(f(x),f(y))_W=(x,y)_V$.
\end{defin}    
Any isometry preserves the norm in the sense that $\|f(x)\|=\|x\|$. Assuming that $x\in {\rm Ker} f$ we obtain $\|x\|=0\, \Rightarrow x=0$ which means that ${\rm Ker} f=\{0\}$ and, therefore,  $f$ is injective. Thus any isometry is injective. The composition of isometries is an isometry.
\begin{defin}
A bijective isometry $f:V\to W$ is called an isomorphism. The unitary spaces related by an isomorphism are called isomorphic and are denoted as $V\cong W$. 
\end{defin} 
We have seen that the Riesz map  $\phi_R: V\to V^*$ is a special isomorphism which eliminates the distinction between covariant and contravariant indices. Therefore, the spaces of tensors  with the same number of indices are isomorphic among themselves, 
$T^n(V)\cong T^{n-1}_1(V)...\cong T^{n-k}_k(V)...\cong T_n(V)$.  Then, in current calculations the position of indices is indifferent and the Einstein convention does not work.

A special case is that of the automorphisms $f\in {\rm Aut}(V)$ that  preserve the inner product of $V$. 
\begin{theor}
The operators $f:V\to V$ that preserve the inner  product as $(f(x),f(y))=(x,y)$ are unitary operators which satisfy $f^{\dagger}=f^{-1}$.
\end{theor}
\begin{demo}
We have $(f(x),f(y))=(f^{\dagger}\circ f(x),y)\, \Rightarrow f^{\dagger}\circ f={\rm id}_V$ and $f$ bijective. Then $f^{\dagger}=f^{-1}$ 
\end{demo} The matrices of the unitary operators in orthonormal bases are unitary, $|f|^+=|f|^{-1}$.
\begin{defin}
The automorphisms that preserve the inner product are called the unitary automorphisms (or unitary transformations). 
\end{defin}
\begin{theor}
The unitary automorphisms $f:V\to V$ transform the orthonormal basis $e$ into the orthonormal basis $e'$ such that $e'_i=f(e_i)$.
\end{theor}
\begin{demo}
First we consider the unitary automorphism $f$ and the basis transformation $e'_i=f(e_i)$.  Then we have $(e'_i,e'_j)=(f(e_i),f(e_j))=(e_i,e_j)=\delta_{ij}$ so that both these bases are orthonormal simultaneously. 
\end{demo} Note that, according to theorem (\ref{ortbasM}), the orthonormal bases are related among themselves only through unitary transformations. 
\begin{defin}
Suppose $V$ is a $n$-dimensional unitary vector space. The unitary automorphisms of $V$ form the unitary group $U(n)\subset {\rm Aut}(V)\cong GL(n,{\Bbb C})$.
\end{defin}

\subsection{Vector spaces with indefinite inner product}

In relativistic physics there are many problems that lead to non-unitary representations of groups and algebras. These can be solved only in complex vector spaces with indefinite metric or indefinite inner product. Moreover, the indefinite metric is related to a specific symmetry which is of high physical interest. 

\subsubsection{Indefinite metric}

Let $V$ be a unitary vector space with its inner product $(~,~)$. In this space we introduce a non-degenerate H-form $H:V\times V\to {\Bbb C}$ as in definition (\ref{H-form}).
\begin{defin}
We say that the operator $h\in {\rm Aut}(V)$ is the metric operator associated to the H-form $H$ if
\begin{equation}\label{funch}
(x,h(y))=H(x,y)\,, \quad \forall\, x,\,y\in V.
\end{equation} 
\end{defin}
\begin{lem}
The metric operator associated to a H-form is selfadjoint, $h^{\dagger}=h$.
\end{lem}
\begin{demo}
$\overline{H(x,y)}=\overline{(x,h(y))}\, \Rightarrow H(y,x)=(h(y),x)=(y,h^{\dagger}(x))$ so that
$(y,h(x))=(y,h^{\dagger}(x)) $ for all $x,\,y\in V$
\end{demo} Consequently, according to the theorem (\ref{spectheor}) there exists at least one orthonormal basis in which the matrix of $h$ takes the canonical form  (\ref{canform}).
\begin{defin}\label{compatible}
The H-form $H$ is compatible with the scalar product of $V$ if $(h(x),h(y))=(x,y)$, i. e. $h\circ h={\rm id}_V$.
\end{defin}   
Under such circumstances, the metric operator $h=h^{\dagger}=h^{-1}$ can have only the eigenvalues $\lambda=\pm 1$ and, therefore, there exist an orthonormal basis, say $\epsilon$,  in which its matrix takes the canonical form
\begin{equation}\label{diah}
| h|_{\epsilon} ={\rm diag}(\, \underbrace{1,1,...,1}_{n_+},\underbrace{-1,-1,...,-1}_{n_-}\, )\,,  
\end{equation} 
where $n_++n_- = {\rm dim}(V)$.  O course, we speak about orthogonality with respect to the inner product which gives the matrix elements 
\begin{equation}\label{hHee}
h_{ij}=(\epsilon_i,h(\epsilon_j))=H(\epsilon_i,\epsilon_j)=\pm \delta_{ij}
\end{equation}
of the matrix $|h|$ in this basis. However, changing this basis into another orthonormal one the canonical form of $|h|$ could change. For this reason we must distinguish among the genuine orthonormal bases with respect to the inner product and the bases in which $|h|$ is diagonal as we show in the next section. 

For the present we remain in the basis $\epsilon$ observing that the metric operator $h$ has only two subspaces of eigenvectors, namely $V_+$ for the eigenvalue $+1$ and $V_-$ for the the eigenvalue $-1$. These subspaces are orthogonal so that $V=V_+\oplus V_-$. In our notations $n_+={\rm dim} V_+$ and  $n_-={\rm dim} V_-$. We exploit this conjecture introducing the  projection operators $p_+$ and $p_-$ such that   $p_+(V)=V_+$ and $p_-(V)= V_-$. These operators are orthogonal, $p_+\circ p_-=0$,   $p_-\circ p_+=0$, form a complete system of projection operators, $p_++p_-={\rm id}_V$ (see definition (\ref{ortproj})), and allow us to write $h=p_+-p_-$. Then we get the general formula 
\begin{equation}
H(x,y)=(x,h(y))=(x,p_+(y))-(x,p_-(y))\,,
\end{equation}
which is basis-independent.

This result is general despite of the apparent restriction imposed by definition (\ref{compatible}). One can convince that considering another approach which starts only with the H-form and has the following steps: definition of the H-form $\Rightarrow$ demonstration of the canonical form (\ref{diah})   $\Rightarrow$ definition of the projector operators $\Rightarrow$ definition of the {\em compatible} genuine inner product $(x,y)=H(x,p_+(y))-H(x,p_-(y))$ of the unitary vector space $V$. Anyway, both these methods are equivalent and leads to the same conclusion that $V_+$ and $V_-$ are orthogonal unitary subspaces. 
\begin{defin}
The numbers $(n_+,n_-)$ gives the {signature} of the  H-form which is its principal characteristic. For this reason we say that two H-forms are equivalent if these have the same signature.
\end{defin}
When the H-form  has an  arbitrary signature, $(n_+,n_-)$,  we say that the metric operator is {\em indefinite} since the quadratic forms 
\begin{equation} 
H(x,x)=|x^1|^2+|x^2|^2+...|x^{n_+}|^2-|x^{n_+ +1}|^2 -....-|x^{n}|^2\,\in{\Bbb R}\,,
\end{equation}
can take any real value. But when the signature is positive, with $n_+ = n$ and $n_- =0$, then $|h|={\bf 1}_{n}$ (since $h={\rm id}_V$) and the space $V$ is a genuine unitary space. Otherwise we must consider explicitly the metric operator $h$.
\begin{defin}
The complex vector space with a metric operator $h$ is called a pseudo-unitary vector space or a vector space with indefinite metric and is denoted by $(V,h)$. Then the H-form $H(~,~)=(~,h(~))$ is called the indefinite inner product. 
\end{defin}
\begin{rem}
The unitary spaces (with inner products) and the pseudo-unitary ones (with indefinite inner products) are called often  inner product spaces and respectively indefinite inner product spaces. 
\end{rem}
The conclusion is that for investigating finite-dimensional vector spaces with indefinite metrics  we have two important tools, the inner products and the H-forms. From the mathematical point of view, the inner product is important for analysing orthogonality properties in the orthonormal bases defined with its help. However, physically  speaking the H-forms giving indefinite inner products are the only quantities getting physical meaning in relativity.

\subsubsection{Dirac conjugation}

Summarising the above results we see that  now we have three different forms: the dual form $\langle~,~\rangle : V^*\times V\to {\Bbb C}$, the inner product $(~,~):V\times V\to {\Bbb C}$ and the H-forms  $H(~,~) : V\times V\to {\Bbb C}$. The first two forms are related through the Riesz representation theorem  while last two are related by Eq. (\ref{funch}).
\begin{theor}
Given a non-degenerate H-form $H$ there exists a bijective anti-linear map $\phi_D : V\to V^*$, called the Dirac isomorphism, such that 
\begin{equation}\label{psih}
\phi_D = \phi_R \circ h\, \Rightarrow \langle \phi_D(x),y\rangle=H(x,y)
\end{equation}
\end{theor}
\begin{demo}
According to the Riesz theorem (\ref{Riesz}) the map $\phi_R$ is bijective and $H$ is non-degenerate. Therefore, $\phi_D$ is bijective.  Then we calculate
\begin{equation}
\langle \phi_D(x),y\rangle=\langle \phi_R(h(x)),y\rangle=(h(x),y)=(x,h(y))=H(x,y)
\end{equation}
\end{demo}
The map $\phi_D$  is uniquely determined by the non-degenerate form $H$ so that we  can say that $\phi_D$ is associated to $H$. 
\begin{defin}
The  H-form  $\overline{H}: V^*\times V^*\to {\Bbb C}$ defined as  
\begin{equation}\label{overHH}
\overline{H}(\hat x,\hat y)=H(\phi_D^{-1}(\hat y),\phi_D^{-1}(\hat x))=\langle \hat y,\phi_D^{-1}(\hat x)\rangle\,.
\end{equation}
is the Dirac adjoint of the H-form $H$. The metric operator $\overline{h}:V^*\to V^*$ that satisfies
\begin{equation}\label{Hhxy}
\overline{H}(\hat x,\hat y)=(\hat x, \overline{h}(\hat y))_*\,,
\end{equation} 
is the Dirac adjoint of the metric operator $h$ of $V$
\end{defin}
The adjoint H-form $\overline{H}$ is  non-degenerate as long as $H$ is non-degenerate and,   $\phi_D$ is bijective.  Moreover, we see that this is correctly defined since is  anti-linear in its first argument, $\overline{H}(\alpha\hat x,\hat y)=\overline{\alpha}\,\overline{H}(\hat x,\hat y)$.

\begin{defin}\label{Dirconj}
The map $\phi_D$ defines the Dirac conjugation with respect to the H-form $H$. We say that  $\overline{x}=\phi_D(x)\in V^*$ is  the Dirac adjoint of the vector $x\in V$ and by $\overline{\hat x}=\phi^{-1}_D(\hat x)$ is the Dirac adjoint of the covector $\hat x\in V^*$.
\end{defin}
\begin{theor}
The Dirac adjoint of the metric operator $h$ is related to its dual $h^*$ as 
\begin{equation}\label{hhhstar}
\overline{h}=h^*=\phi_R\circ h \circ \phi_R^{-1}
\end{equation}
\end{theor} 
\begin{demo}
We start with Eq. (\ref{overHH}), $\overline{H}(\hat x,\hat y)=H(\phi_D^{-1}(\hat y),\phi_D^{-1}(\hat x))=(\phi_D^{-1}(\hat y),h\circ\phi_D^{-1}(\hat x))=(h\circ  \phi_R^{-1}(\hat y),\phi_R^{-1}(\hat x))$ while from Eq. (\ref{Hhxy}) we have $\overline{H}(\hat x,\hat y)=(\hat x, \overline{h}(\hat y))_*=(\phi_R^{-1}\circ \overline{h}(\hat y), \phi_R^{-1}(\hat x))$. We compare these results and use Eq. (\ref{ffstarcr}) obtaining the desired formula
\end{demo} One can read this result on the diagram:
\begin{equation}\label{dia7}
\begin{array}{rcl}
&h=h^{\dagger}&\\
V~&\overrightarrow{\hspace*{20mm}}&V\\
\left.\begin{array}{c}
\\
\phi_R\\
\\
\end{array} \right\downarrow~&&
\left\downarrow
\begin{array}{c}
\\
\phi_R\\
\\
\end{array}\right.\\
V^*&\overrightarrow{\hspace*{20mm}}&V^*\\
&h^*=\overline{h}&
\end{array}
\end{equation}
\begin{rem} The conclusion is that if  $(V,h)$ is a vector space with indefinite metric $h$  then its dual $(V^*,h^*)$ is of indefinite metric too having the metric operator $h^*$.
\end{rem}
We specify that these operators have the same matrix since $|h^*|=|h|$ (see definition (\ref{matf})). Another consequence of Eq. (\ref{hhhstar}) is  that $\phi_D^*=\phi_D$ as we can prove calculating $\phi_D^*=h^*\circ \phi_R=\phi_R\circ h$ (since $\phi_R^*=\phi_R$). However, this expected result verifies the consistency of our approach where all  the mappings of ${\rm Lin}(V,V^*)$ must be selfdual as we stated in definition (\ref{selfdual}).
\begin{defin}
We say that  $\overline{f}:V\to V$ is the Dirac adjoint of the operator $f:V\to V$if 
\begin{equation}\label{HHffhh}
H(x,f(y))=H(\overline{f}(x),y)\, \Rightarrow \, \overline{f}=h\circ f^{\dagger}\circ h\,. 
\end{equation}
\end{defin}
The Dirac conjugation has the obvious rules: $\overline{(\overline{f})}=f$ and
\begin{equation}
\overline{(f+g)}=\overline{f}+\overline{g}\,,\quad\overline{ (f\circ g)}=\overline{g}\circ \overline{f}\,, \quad  (\alpha f)^{\dagger}=\overline{\alpha}\overline{f}\,,\,\,\alpha\in {\Bbb C}\,.
\end{equation}
It follows naturally to give the definition:
\begin{defin}
We say that the Dirac self-adjoint (or simply self-adjoint) operators satisfy $f=\overline{f}$. 
\end{defin}
We observe that the Dirac self-adjoint operators $f$ can be related to the operators $f\circ h$ which are selfadjoint with respect to the inner product $(~,~)$ of $V$. Then we may exploit the fact that these last mentioned operators have spectral representations.
\begin{theor}
Any Dirac selfadjoint operator $f$ has the spectral representation
\begin{equation}
f=\sum_{\lambda\in {\rm Sp}(f\circ h)}\lambda \,p_{\lambda}\circ h\,.
\end{equation}
\end{theor}  
\begin{demo}
The selfadjoint operator $f\circ h$ has real eigenvalues and the eigenspaces $V_{\lambda}$ whose projection operators give the spectral representation  (\ref{spreprep}). Then we use $f=(f\circ h)\circ h$
\end{demo}

Finally let us see which is the role of the pseudo-unitary operators.
\begin{theor}\label{AutVh}
The automorphisms $f\in{\rm Aut}(V,h)$ that preserve the indefinite metric $h$ of signature $(n_+,n_-)$ are pseudo-unitary operators that accomplish $f^{-1}=\overline{f}$. These automorphisms form the pseudo-unitary group $U(n_+,n_-)$.
\end{theor}
\begin{demo}
The transformed H-form is $H(f(x),f(y))=H(x,\overline{f}\circ f(y))=H(x,y)$ for all $x.y\in V$ 
$\Rightarrow \overline{f}\circ f={\rm id}_V$. 
\end{demo} In other words, we found the isomorphism ${\rm Aut}(V,h)\cong U(n_+,n_-)$.

\subsubsection{Representations in $h$-orthonormal dual bases} 

Our purpose is to represent the pseudo-unitary space $(V,h)$ as an arithmetic space  ${\Bbb C}^n$ with $n={\rm dim} V$. We start defining the bases we need.
\begin{defin}
We say that a pair of dual bases $(e,\hat e)$ are $h$-orthonormal dual bases if 
\begin{equation}\label{psih1}
\phi_D(e_i)=h_{ij}\hat e^j\,,  \quad \phi^{-1}_D(\hat e^i)=h^{ij} e_j\,, 
\end{equation}
where we denote $h_{ij}=h^{ij}=\pm\delta_{ij}$ the canonical form given by Eq, (\ref{hHee}) in order to preserve the Einstein convention. 
\end{defin}
In general  the $h$-orthonormal bases are different from the  orthonormal ones (with respect to the inner product). However, the orthonormal basis in which the H-form is brought in diagonal form is in the same time $h$-orthonormal.  Moreover, these two type of  bases can coincide  as sub-bases inside  the orthogonal unitary subspaces $V_+$ and $V_-$.  
\begin{rem}
The conclusion is that there is a family of  $h$-orthonormal bases which are orthonormal with respect to the scalar product.
\end{rem}
Obviously, there is another family of $h$-orthonormal bases which are no longer orthonormal. We shall see later how these properties are related to some subgroups of the  isometry groups.

From Eqs. (\ref{psih1}) we deduce that in these bases we have the identities
\begin{equation}
H(e_i,e_j)=\overline{H}(\hat e^i,\hat e^j)=h^{ij}=h_{ij}\,.
\end{equation}
Then we see that the standard expansions $x=x^ie_i$ and $\hat y =\hat y_i\hat e ^i$ give the expansions of the Dirac conjugated vectors and covectors, 
\begin{equation}\label{calc1}
\overline{x}=\phi_D(x)=\overline{x^i}\phi_D(e_i)=h_{ij}\overline{x^i}\hat e^j\,, \quad \overline{\hat y}=\phi_D^{-1}(\hat y)= \overline{\hat y_i}\phi_D^{-1} (\hat e^i)= {h}^{\,ij}\overline{\hat y_i}e_j\,,
\end{equation}
as resulted from Eqs.  (\ref{psih1}) and definition (\ref{Dirconj}). Hereby, we deduce   
\begin{equation}
\langle \overline{x}, y\rangle=H(x,y)=h_{ij}\overline{x^i}y^j\,,\quad
\langle \hat{x}, \overline{\hat y}\rangle=\overline{H}(\hat y,\hat x)={h}^{ij}\hat x_i \overline{\hat y_j}\,.
\end{equation}
These expansions lay out the components of the vectors and covectors defined canonically as 
\begin{eqnarray}
x^i &=&\langle \hat e^i,x\rangle=h^{ij}\langle \phi_D(e_j),x)=h^{ij}H(e_j,x)\,,\\
\hat y_i &=& \langle \hat y, e_i\rangle=h_{ij}\langle \hat y,\phi_D^{-1}\,.(e_j)\rangle=h_{ij}\overline{H}(\hat e^j,\hat y)
\end{eqnarray}
These equations suggest the opportunity of introducing covariant component  for vectors,  $x_i=H(e_i,x)=h_{ij}x^j$, and contravariant ones for covectors $\hat y^i=\overline{H}(\hat e^i,\hat y)=h^{ij}\hat y_j$.  Moreover, we observe that  the  matrix elements of any operator can be written as $f^i_j=h^{ik}f_{kj}$ where
\begin{equation}
f_{ij}=H( e_i, f(e_j))=h_{ik}\langle \hat e^k, f(e_j)\rangle=h_{ik}f_{j}^k\,,
\end{equation}  
since our excessive notation ($h_{ij}=h^{ij}$) helps us to preserve the Einstein convention. In fact,  the metric operator can raise or lower indices as for example in the mapping 
\begin{equation}
\phi_D\otimes...\otimes \phi_D^{-1}.... : V\otimes V.....\otimes V^*\otimes V^*....\to
V^*\otimes V.....\otimes V\otimes V^*...
\end{equation}
which lowers the first index and rises another one giving the new components
\begin{equation}
 \tau^{..i_2..l_2}_{k_1..j_1..}=h_{i_1k_1}...h^{j_2l_2}\cdots \tau^{i_1 i_2...}_{.....j_1j_2..}
\end{equation} 
This means that the Dirac isomorphism and its inverse map are able to relate among themselves all the tensors with the same number of components, $T^n(V)\cong ...T^{n-k}_k(V)...\cong T_n(V)$. 
\begin{defin}
A physical tensor quantity of rank $n$  is represented by all the tensors with $n$ indices, indifferent of their positions.   
\end{defin}  

The vectors,  covectors and their Dirac conjugated are represented by matrices as
\begin{eqnarray}
&&x\in V \to |x\rangle = {\rm Rep}_e (x)\,, \quad \hat y\in V^*\to \langle \hat y|={\rm \hat Rep}_{\hat e} (\hat y)\,,\\
&&\overline{x}\in V^*\to \langle \overline{x}| ={\rm \hat Rep}_{\hat e} (\overline{x})\,, \quad 
\overline{\hat y}\in V\to | \overline{\hat y}\rangle ={\rm Rep}_e (\overline{\hat y})\,,
\end{eqnarray} 
\begin{defin}
The Dirac conjugation of the representation matrices of the vectors and covectors read 
\begin{equation}
\overline{{\rm Rep}_e (x)}={\rm \hat Rep}_{\hat e} (\overline{x})\,,\quad \overline{{\rm \hat Rep}_{\hat e} (\hat y)}={\rm Rep}_{e} (\overline{\hat y})\,.
\end{equation}
\end{defin}
\begin{theor}
In any matrix representation of a pair of $h$-orthonormal bases $(e,\hat e)$ where $|i\rangle ={\rm Rep}_{e}(e_i)$ and $\langle i | ={\rm \hat Rep}_{\hat e}(\hat e^i)$ the basic calculation  rules of the Dirac conjugation are given by $\overline{ | i\rangle}=\langle i | \,|h|$ and $ \overline{\langle i|}=|h|\,|i\rangle$. 
\end{theor} 
\begin{demo}
We use Eq. (\ref{calc1}) and the above definition calculating 
\begin{equation}
\overline{|i\rangle}=\overline{{\rm Rep}_e (e_i)}= {\rm \hat Rep}_{\hat e} (\overline{e_i})=\overline{h_{ij}}\,
{\rm \hat Rep}_e (\hat e^j)=\langle i|\,|h| 
\end{equation}  
since in these bases $|h|$ is a diagonal real-valued  matrix. The second rule has to be proved similarly.
\end{demo}
\begin{rem}
In concrete calculations we can exploit the {natural} properties of  the ket and bra bases  of this representation, $(| i\rangle)^+=\langle i |$ and $(\langle i|)^+=|i\rangle$. 
\end{rem}
For example, we can write
\begin{equation}
\overline{|x\rangle}=\langle \overline{x}| =(|x\rangle)^+|h|\,, \quad 
\overline{\langle \hat y|}= | \overline{\hat y}\rangle =|h|(\langle \hat y |)^+\,.
\end{equation}
With these notations, the indefinite inner product can be represented as a product of  matrices\,, 
\begin{equation}
H(x,y)=\langle \overline{x}, y\rangle =\langle \overline{x}| y\rangle= (\langle\overline{x}|)(| y\rangle)=(|x\rangle)^+|h| \,|y\rangle=\langle x|h|y\rangle\,,
\end{equation}
according to the Dirac's "bra-ket" notation presented in section 2.2.2.. These simple rules lead to the wrong idea that the Dirac conjugation is an occasional trick which  reduces to the calculation of the Hermitian adjoint multiplied with a certain diagonal matrix.  In fact, we have seen that this is the most general definition of a conjugation in complex vector spaces that includes the Hermitian conjugation as a particular case when $|h|={\bf 1}$. 

The linear operators $f\in {\rm End} (V)$ are represented by their matrices in  $h$-orthonormal bases,  $|f|={\rm Rep}_{(e,e)}(f)\in {\rm End} ({\Bbb C}^n)$ ($n={\rm dim}(V)$), whose matrix elements are given by Eq. (\ref{fffReee}). 
Then we have the representation  
\begin{equation}
H(x,f(y))=\langle \overline{x}|{\rm Rep}_{(e,e)}(f)|y\rangle=\langle \overline{x}|f|y\rangle=\langle x|hf|y\rangle\,.
\end{equation}
In a pair of $h$-orthonormal bases we can write spectral representations. Supposing that the signature of $h$ is $(n_+,n_-)$ we can write the matrices representing of the projection operators 
\begin{equation}
|p_+|={\rm Rep}_{(e,e)}(p_+)=\sum_{k=1}^{n_+}|k\rangle\langle k| \,,\quad 
|p_-|={\rm Rep}_{(e,e)}(p_-)=\sum_{k=n_++1}^{n}|k\rangle\langle k| \,,
\end{equation}   
giving the representation of the metric operator
\begin{equation}
{\rm Rep}_{(e,e)}(h)=|h|=|p_+|-|p_-|\,.   
\end{equation}
In general, a linear operator $f\in {\rm End}(V)$ has the  representation
\begin{equation}
|f|={\rm Rep}_{(e,e)}(f)=\sum_{i=1}^{n}\sum_{j=1}^{n}|i\rangle f_{j}^i\langle j|\,.
\end{equation}

The Hermitian properties of these matrices are determined (indirectly) by the H-form that defines the Dirac conjugation. 
\begin{theor}
The matrices in $h$-orthonormal bases of two Dirac conjugated operators, $f$ and $\overline{f}$, are related as 
\begin{equation}
{\rm Rep}_{(e,e)}(\overline{f})= |h| \,{\rm Rep}_{(e,e)}(f)^+|h|\,.
\end{equation}  
\end{theor}
\begin{demo}
We represent Eq. (\ref{HHffhh})
\end{demo}
This property justifies the common notation ${\rm Rep}_{(e,e)}(\overline{f})=\overline{{\rm Rep}_{(e,e)}(f)}=\overline{|f|}$ that allows us to write simply, $\overline{|f|}=|h|\,|f|^+|h|$. In this representation, the matrices of the Dirac selfadjoint operators satisfy  $|f|= |h| \,|f|^+|h|$ while for the pseudo-unitary operators  we have $|f|^{-1}= |h| \,|f|^+|h|$.

We have shown that the pseudo-unitary operators form the group of automorphisms ${\rm Aut}(V,h)\cong U(n_+,n_-)$ (see theorem (\ref{AutVh})).
\begin{theor}
The automorphisms of $U(n_+,n_-)$ transform any $h$-orthonormal basis into a $h$-orthonormal basis.
\end{theor}
\begin{demo}
Let $e'_i=f(e_i)$ the vectors of the new basis. Then we have $H(e'_i,e'_j)=H(f(e_i),f(e_j))=H(e_i,e_j)$ since  $f$ is pseudo-unitary
\end{demo}
\begin{defin}
The matrices $|f|={\rm Rep}_{(e,e)}(f)$ of the pseudo-unitary operators form the fundamental representation of the group $U(n_+,n_-)$.
\end{defin}
We shall see that the fundamental representation is in fact the representation in which this group is defined.

\subsubsection{Euclidean and pseudo-Euclidean spaces}

The real vector spaces have similar properties as the complex ones  but with  real numbers taking over the role of the complex ones. Thus the inner product becomes {\em symmetric} and  linear in both its terms such that the adjoint operator of the operator $f$ is the {\em transposed} operator $f^T$ which satisfies $(f^T(x),y)=(x,f(y))$. The H-forms become now  the symmetric forms $\eta:V\times V\to {\Bbb R}$  associated to the symmetric metric operators denoted with the same symbol, $\eta=\eta^T=\eta^{-1}$, and defined as
\begin{equation}
\eta(x,y)=(x,\eta(y))\,.
\end{equation} 
This metric operator can be brought in diagonal canonical form of a given signature $(n_+,n_-)$ with $n_++n_-=n={\rm dim} V$.  
\begin{defin}
The real vector spaces with inner products are Euclidean spaces while those with indefinite metrics are called pseudo-Euclidean. We denote the pseudo-Euclidean space by $(V,\eta)$ where $\eta$ is the metric of  signature $(n_+,n_-)$.    
\end{defin}
Any space with indefinite metric has a compatible inner product defining the orthogonal properties. We must specify that now the Riesz and Dirac isomorphisms are {\em linear} so that the Riesz map can be confused with the natural duality which is independent on the inner product or indefinite metric we use.

\begin{defin}
A pair of dual bases $(e,\hat e)$ are $\eta$-orthonormal dual bases if 
\begin{equation}\label{psih1}
\phi_D(e_i)=\eta_{ij}\hat e^j\,,  \quad \phi^{-1}_D(\hat e^i)=\eta^{ij} e_j\,, 
\end{equation}
where now $\eta_{ij}=\eta^{ij}=\pm\delta_{ij}$ the canonical form similar to those given by Eq, (\ref{hHee}). 
\end{defin}
As in the case of pseudo-unitary spaces, there are $\eta$-orthonormal bases that can be sometimes   simultaneously  orthonormal with respect to the inner product. 

The standard expansions $x=x^ie_i$ and $\hat y =\hat y_i\hat e ^i$ give the expansions of the Dirac conjugated vectors and covectors, 
\begin{equation}\label{calc1}
\overline{x}=\phi_D(x)={x^i}\phi_D(e_i)=\eta_{ij}{x^i}\hat e^j\,, \quad \overline{\hat y}=\phi_D^{-1}(\hat y)= {\hat y_i}\phi_D^{-1} (\hat e^i)= {\eta}^{\,ij}{\hat y_i}e_j\,,
\end{equation}
since all the components are real-valued.  Hereby, we deduce   
\begin{equation}
\langle \overline{x}, y\rangle=\eta(x,y)=\eta_{ij}\overline{x^i}y^j\,,\quad
\langle \hat{x}, \overline{\hat y}\rangle={\eta}(\hat y,\hat x)={\eta}^{ij}\hat x_i {\hat y_j}\,,
\end{equation}
and we recover the mechanism of raising or lowering indices presented in the previous section but now using the metric $\eta$.

For each pair of $\eta$-orthonormal bases we obtain a representation in the space ${\Bbb R}^n$ with the same rules as in the case of the complex vector spaces. The difference is that now the Dirac adjoint read
\begin{equation}
\overline{|x\rangle}=\langle \overline{x}| =(|x\rangle)^T|\eta|\,, \quad 
\overline{\langle \hat y|}= | \overline{\hat y}\rangle =|\eta|(\langle \hat y |)^T\,,
\end{equation}
and  the indefinite inner product can be represented as a product of  matrices and finally written in the "bra-ket" notation as, 
\begin{equation}
\eta(x,y)=\langle \overline{x}, y\rangle =\langle \overline{x}| y\rangle= (\langle\overline{x}|)(| y\rangle)=(|x\rangle)^T|\eta| \,|y\rangle=\langle x|\eta|y\rangle\,.
\end{equation}

The operators of ${\rm End}(V,\eta)$ have similar properties with those defined on unitary or pseudo-unitary spaces tking into account that now the Dirac adjoint of the operator $f$ reads $\overline{f}=\eta \circ f^T\circ \eta$ while the terminology is changed as\\ 

\noindent\begin{tabular}{lcl}
Complex inner product spaces & & Real inner product spaces\\
&&\\
selfadjoint,  $f=f^{\dagger}$ &$\Rightarrow$ & selfadjoint, $f=f^T$\\
Dirac selfadjoint, $f=h\circ f^{\dagger}\circ h$&$\Rightarrow$ &(Dirac) selfadjoint, $f=\eta\circ f^T\circ \eta$\\
unitary, $f^{-1}=f^{\dagger}$&$\Rightarrow$ & {\em orthogonal}, $f^{-1}=f^T$\\
pseudo-unitary, $f^{-1}=h\circ f^{\dagger}\circ h$&$\Rightarrow$ &{\em pseudo-orthogonal}, $f^{-1}=\eta\circ f^T\circ\eta$
\end{tabular}\\

The linear operators $f\in {\rm End} (V,\eta)$ are represented by their matrices in  $\eta$-orthonormal bases,  $|f|={\rm Rep}_{(e,e)}(f)\in {\rm End} ({\Bbb R}^n)$ ($n={\rm dim}(V)$), whose matrix elements are given by Eq. (\ref{fffReee}). 
Then we can write  
\begin{equation}
\eta(x,f(y))=\langle \overline{x}|{\rm Rep}_{(e,e)}(f)|y\rangle=\langle \overline{x}|f|y\rangle=\langle {x}|\eta f|y\rangle\,.
\end{equation}
Moreover, we recover the spectral representation in a pair of $\eta$-orthonormal bases as in the case of the complex vector spaces. 

The orthogonal properties of the representation matrices are determined  by the metric $\eta$ that defines the Dirac conjugation. 
\begin{theor}
The matrices in $\eta$-orthonormal bases of two Dirac conjugated operators, $f$ and $\overline{f}$, are related as 
\begin{equation}
{\rm Rep}_{(e,e)}(\overline{f})= |\eta| \,{\rm Rep}_{(e,e)}(f)^T|\eta|\,.
\end{equation}  
\end{theor}
This property enables us to denote ${\rm Rep}_{(e,e)}(\overline{f})=\overline{{\rm Rep}_{(e,e)}(f)}=\overline{|f|}$ such that $\overline{|f|}=|\eta|\,|f|^T|\eta|$. In this representation, the matrices of the Dirac selfadjoint operators satisfy  $|f|= |\eta| \,|f|^T|\eta|$ while for the pseudo-orthogonal operators  we have $|f|^{-1}= |\eta| \,|f|^T|\eta|$.

\begin{theor}
If $(V,\eta)$ is a  pseudo-Euclidean space and $\eta$ has the signature $(n_+, n_-)$ then ${\rm Aut}(V,h)$ is isomorphic to the pseudo-orthogonal group $O(n_+, n_-)\subset {\rm Aut}(V)\cong GL(n,{\Bbb R})$.
\end{theor}  
\begin{demo}
This is a particular case of theorem (\ref{AutVh}) restricted to pseudo-orthogonal automorphisms 
\end{demo}
The pseudo-orthogonal automorphisms preserve the metric and transform any $\eta$-orthonormal basis into another basis of the same type.

\end{document}